\documentclass[11pt]{article}
\usepackage{amssymb}
\usepackage[dvips]{epsfig}
\usepackage{times}
\usepackage{graphicx}

\newtheorem{Lemma}{Lemma}[section]
\newtheorem{Theorem}{Theorem}

\newtheorem{Remark}{Remark}
\newtheorem{Assumption}{Assumption}

\newenvironment{Proof}%
 {\begin{trivlist} \item[]{\bf Proof. }}%
 {\hspace*{\fill}$\rule{.3\baselineskip}{.35\baselineskip}$\end{trivlist}}
 {\begin{trivlist}\item[]\textbf{Acknowledgments }}{\end{trivlist}}

\makeatletter \@addtoreset{equation}{section} \makeatother

\newcommand{\R}{\mathbb{R}}
\newcommand{\Z}{\mathbb{Z}}

\newfam\bifam
\font\tenbi=cmmib10 scaled \magstep1 \font\sevenbi=cmmib10 at 11pt
\font\fivebi=cmmib10 at 6pt \textfont\bifam = \tenbi
\scriptfont\bifam= \sevenbi \scriptscriptfont\bifam= \fivebi

\begin{document}

\title{\bf Normal form for travelling kinks \\ in discrete Klein--Gordon lattices}

\author{Gerard Iooss$^\dagger$ and Dmitry E. Pelinovsky$^{\dagger \dagger}$\\
{\small $^{\dagger}$ Institut Universitaire de France, INLN, UMR
CNRS-UNSA 6618, 1361 route des Lucioles, 06560 Valbonne, France}  \\
{\small $^{\dagger\dagger}$ Department of Mathematics, McMaster
University, 1280 Main Street West, Hamilton, Ontario, Canada, L8S
4K1} }

\date{\today}
\maketitle

\begin{abstract}
We study travelling kinks in the spatial discretizations of the 
nonlinear Klein--Gordon equation, which include the discrete 
$\phi^4$ lattice and the discrete sine--Gordon lattice. The 
differential advance-delay equation for travelling kinks is reduced 
to the normal form, a scalar fourth-order differential equation, 
near the quadruple zero eigenvalue. We show numerically 
non-existence of monotonic kinks (heteroclinic orbits between 
adjacent equilibrium points) in the fourth-order equation. Making 
generic assumptions on the reduced fourth-order equation, we prove 
the persistence of bounded solutions (heteroclinic connections 
between periodic solutions near adjacent equilibrium points) in the 
full differential advanced-delay equation with the technique of 
center manifold reduction. Existence and persistence of multiple 
kinks in the discrete sine--Gordon equation are discussed in 
connection to recent numerical results of \cite{ACR03} and results 
of our normal form analysis. 
\end{abstract}

\section{Introduction}

Spatial discretizations of the nonlinear partial differential
equations represent discrete dynamical systems, which are
equivalent to chains of coupled nonlinear oscillators or discrete
nonlinear lattices. Motivated by various physical applications and
recent advances in mathematical analysis of discrete lattices, we
consider the discrete Klein--Gordon equation in the form:
\begin{equation}
\label{Klein-Gordon} \ddot{u}_n = \frac{u_{n+1} - 2 u_n +
u_{n-1}}{h^2} + f(u_{n-1},u_n,u_{n+1}),
\end{equation}
where $u_n(t) \in \R$, $n \in \Z$, $t \in \R$, $h$ is the lattice
step size, and $f(u_{n-1},u_n,u_{n+1})$ is the nonlinearity
function. The discrete lattice (\ref{Klein-Gordon}) is a
discretization of the continuous Klein--Gordon equation, which
emerges in the singular limit $h \to 0$:
\begin{equation}
\label{KG} u_{tt} = u_{xx} + F(u),
\end{equation}
where $u(x,t) \in \R$, $x \in \R$, and $t \in \R$. In particular,
we study two versions of the Klein--Gordon equation (\ref{KG}),
namely the $\phi^4$ model
\begin{equation}
\label{phi4} u_{tt} = u_{xx} + u (1 - u^2)
\end{equation}
and the sine--Gordon equation
\begin{equation}
\label{sine-Gordon} u_{tt} = u_{xx} + \sin(u).
\end{equation}
We assume that the spatial discretization $f(u_{n-1},u_n,u_{n+1})$
of the nonlinearity function $F(u)$ is symmetric,
\begin{equation}
\label{nonlinearity-symmetry} f(u_{n-1},u_n,u_{n+1}) =
f(u_{n+1},u_n,u_{n-1}),
\end{equation}
and consistent with the continuous limit,
\begin{equation}
\label{nonlinearity-continuous} f(u,u,u) = F(u).
\end{equation}
We assume that the zero equilibrium state always exists with $F(0)
= 0$ and $F'(0) = 1$. This normalization allows us to represent
the nonlinearity function in the form:
\begin{equation}
\label{nonlinearity-representation} f(u_{n-1},u_n,u_{n+1}) = u_n +
Q(u_{n-1},u_n,u_{n+1}),
\end{equation}
where the linear part is uniquely normalized (parameter $h$ can be 
chosen so that the linear term of $u_{n-1} + u_{n+1}$ is cancelled) 
and the nonlinear part is represented by the function 
$Q(u_{n-1},u_n,u_{n+1})$. In addition, we assume that (i) $F(u)$ and 
$Q(u_{n-1},u_n,u_{n+1})$ are odd such that $F(-u) = -F(u)$ and 
$Q(-v,-u,-w) = - Q(v,u,w)$ and (ii) a pair of non-zero equilibrium 
points $u_+ = - u_- \neq 0$ exists, such that
\begin{equation}
\label{nonlinearity-kink} F(u_{\pm}) = 0, \qquad F'(u_{\pm}) < 0,
\end{equation}
while no other equilibrium points exist in the interval $u \in
[u_-,u_+]$. For instance, this assumption is verified for the
$\phi^4$ model (\ref{phi4}) with $u_{\pm} = \pm 1$ and for the
sine--Gordon equation (\ref{sine-Gordon}) with $u_{\pm} = \pm \pi$. 

We address the fundamental question of existence of travelling wave 
solutions in the discrete Klein--Gordon lattice 
(\ref{Klein-Gordon}). Since discrete equations have no translational 
and Lorentz invariance, unlike the continuous Klein--Gordon equation 
(\ref{KG}), existence of travelling waves, pulsating travelling 
waves and travelling breathers represents a challenging problem of 
applied mathematics (see recent reviews in \cite{S03,IJ04}). 
        
Our work deals with the travelling kinks between the non-zero 
equilibrium states $u_{\pm}$. We are not interested in the 
travelling breathers and pulsating waves near the zero equilibrium 
state since the zero state is linearly unstable in the dynamics of 
the discrete Klein--Gordon lattice (\ref{Klein-Gordon}) with $F'(0) 
= 1 > 0$. Indeed, looking for solutions in the form $u_n(t) = e^{i 
\kappa h n + \lambda n}$, we derive the dispersion relation for 
linear waves near the zero equilibrium state:
$$
\lambda^2 = 1 - \frac{4}{h^2} \sin^2 \frac{\kappa h}{2}.
$$
It follows from the dispersion relation that there exists 
$\kappa_*(h) > 0$ such that $\lambda^2 > 0$ for $0 \leq \kappa < 
\kappa_*(h)$. On the other hand, the non-zero equilibrium states 
$u_+$ and $u_-$ are neutrally stable in the dynamics of the discrete 
Klein--Gordon lattice (\ref{Klein-Gordon}) with $F'(u_{\pm}) < 0$. 
We focus hence on bounded heteroclinic orbits which connect the 
stable non-zero equilibrium states $u_-$ and $u_+$ in the form:
\begin{equation}
\label{kink-shape} u_n(t) = \phi(z), \qquad z = hn - ct,
\end{equation}
where the function $\phi(z)$ solves the differential advance-delay
equation:
\begin{equation}
\label{advance-delay} c^2 \phi''(z) = \frac{\phi(z+h)-2 \phi(z) +
\phi(z-h)}{h^2} + \phi(z) + Q(\phi(z-h),\phi(z),\phi(z+h)).
\end{equation}
We consider the following class of solutions of the differential
advance-delay equation (\ref{advance-delay}): (i) $\phi(z)$ is
twice continuously differentiable function on $z \in \mathbb{R}$;
(ii) $\phi(z)$ is monotonically increasing on $z \in \mathbb{R}$
and (iii) $\phi(z)$ satisfies boundary conditions:
\begin{equation}
\label{kink-bc} \lim_{z \to -\infty} \phi(z) = u_-, \qquad \lim_{z
\to +\infty} \phi(z) = u_+.
\end{equation}

It is easy to verify that the continuous Klein--Gordon equation 
(\ref{KG}) with $F(u)$ in (\ref{nonlinearity-kink}) yields a 
travelling kink solution in the form $u = \phi(z)$, $z = x-ct$ for 
$|c| < 1$. However, the travelling kink can be destroyed in the 
discrete Klein--Gordon lattice (\ref{Klein-Gordon}), which results 
in violation of one or more conditions on $\phi(z)$. For instance, 
the bounded twice continuously differentiable solution $\phi(z)$ may 
develop non-vanishing oscillatory tails around the equilibrium 
states $u_{\pm}$ \cite{IJ04}.           

A recent progress on travelling kinks was reported for the discrete 
$\phi^4$ model. Four particular spatial discretizations of the 
nonlinearity $F(u) = u(1-u^2)$ were proposed with four independent 
and alternative methods \cite{BT97,S97,FZK99,K03}, where the 
ultimate goal was to construct a family of translation-invariant 
stationary kinks for $c = 0$, that are given by continuous, 
monotonically increasing functions $\phi(z)$ on $z \in \mathbb{R}$ 
with the boundary conditions (\ref{kink-bc}). Exceptional 
discretizations were generalized in \cite{BOP05,DKY05a,DKY05b}, 
where multi-parameter families of cubic polynomials 
$f(u_{n-1},u_n,u_{n+1})$ were obtained. It was observed in numerical 
simulations of the discrete $\phi^4$ model \cite{S97} that the 
effective translational invariance of stationary kinks implies 
reduction of radiation diverging from moving kinks. New effects such 
as self-acceleration were reported for some of the exceptional 
discretizations in \cite{DKY05a,DKY05b}. Nevertheless, from a 
mathematically strict point of view, we shall ask if exceptional 
discretizations guarantee existence of true travelling kinks 
(heteroclinic orbits) at least for small values of $c$. The question 
was answered negatively in \cite{OPB05}, where numerical analysis of 
beyond-all-order expansions was developed. It was shown that 
bifurcations of travelling kink solutions from $h = 0$ to small 
non-zero $h$ do not generally occur in the discrete $\phi^4$ model 
with small values of $c$, even if the exceptional discretizations 
allow these bifurcations for $c = 0$. It was also discovered in 
\cite{OPB05} that bifurcations of travelling kink solutions may 
occur for finitely many isolated values of $c$ far from the limit $c 
= 0$.

The discrete sine-Gordon model was also subject of recent studies. 
Numerical computations were used to identify oscillatory tails of 
small amplitudes on the travelling kink solutions \cite{EF90,SZE00}. 
These tails were explained with analysis of central manifold 
reductions (carried out without the normal form reductions) 
\cite{FR05}. Exceptional discretizations of the nonlinearity $F(u) = 
\sin(u)$ were suggested by the topological bound method in 
\cite{SW94} (see the review in \cite{S99}) and by the "inverse" 
(direct substitution) method in \cite{FZK99}. Simultaneously with 
the absence of single kinks in the sine-Gordon lattices, multiple 
kinks (between non-zero equilibrium points of $\pi ({\rm mod} 2 
\pi)$) were discovered with the formal reduction of the discrete 
lattice (\ref{Klein-Gordon}) to the fourth-order ODE problem in 
\cite{CK00} and confirmed with numerical analysis of the 
differential advanced-delay equation in \cite{ACR03}.

Our work is motivated by the recent advances in studies of 
differential advanced-delay equations from the point of dynamical 
system methods such as central manifold reductions and normal forms 
(see pioneer works in \cite{I00,IK00}). The same methods were 
recently applied in \cite{PR05} to travelling solitary waves in 
discrete nonlinear Schr\"{o}dinger lattices near the maximum group 
velocity of linear waves. It was shown in \cite{PR05} that 
non-existence of travelling solitary waves can be derived already in 
the truncated (polynomial) normal form. We shall exploit this idea 
to the class of travelling kinks in discrete Klein--Gordon lattices 
near the specific speed $c = 1$ and the continuous limit $h = 0$. 
This particular point corresponds to the quadruple zero eigenvalue 
on the central manifold of the dynamical system. A general 
reversible normal form for the quadruple zero eigenvalue was derived 
and studied in \cite{I95} (see also tutorials in the book 
\cite{IA98}). 

We note that the discrete Klein--Gordon lattice was considered in 
\cite{IK00} but the nonlinearity $F(u)$ was taken to be decreasing 
near $u = 0$, such that the quadruple zero eigenvalue was not 
observed in the list of possible bifurcations of travelling wave 
solutions (see Figure 1 in \cite{IK00}). On the other hand, the 
quadruple zero eigenvalue occurred in the discrete 
Fermi--Pasta--Ulam lattice studied in \cite{I00}, where the symmetry 
with respect to the shift transformation was used to reduce the 
bifurcation problem to the three-dimensional center manifold. Since 
the reversible symmetry operator for the Fermi-Pasta-Ulum lattice is 
minus identity times the reversibility symmetry for the 
Klein--Gordon lattice, the reversible normal forms for quadruple 
zero eigenvalue are different in these two problems. We focus here 
only on the case of the discrete Klein--Gordon lattice.

Our strategy is as follows. We shall develop a decomposition of 
solutions of the differential advance-delay equation into two parts: 
a four-dimensional projection to the subspace of the quadruple zero 
eigenvalue and an infinite-dimensional projection to the hyperbolic 
part of the problem. By using a suitable scaling, we shall truncate 
the resulting system of equations with a scalar fourth-order 
equation, which is similar to the one formally derived in 
\cite{CK00}. We refer to this fourth-order equation as to the normal 
form for travelling kinks. We shall develop a numerical analysis of 
the fourth-order equation and show that no monotonic heteroclinic 
orbits from $u_-$ to $u_+$ exist both in the discrete $\phi^4$ and 
sine-Gordon lattices. Rigorous persistence analysis of bounded 
solutions (heteroclinic orbits between periodic solutions near $u_-$ 
and $u_+$) is developed with the technique of center manifold 
reduction. Our main conclusion is that the differential 
advance-delay equation (\ref{advance-delay}) has no monotonic 
travelling kinks near the point $c = 1$ and $h = 0$ but it admits 
families of non-monotonic travelling kinks with oscillatory tails.

It seems surprising that the truncated normal form is independent on 
the discretizations of the nonlinearity function 
$f(u_{n-1},u_n,u_{n+1})$ and the negative result extends to all 
exceptional discretizations constructed in 
\cite{SW94,S97,S99,BT97,FZK99,K03,BOP05,DKY05a,DKY05b}. Any 
one-parameter curves of the monotonic travelling kinks, which 
bifurcate on the plane $(c,h)$ from the finite set of isolated 
points $(c_*,0)$ (see numerical results in \cite{OPB05}), may only 
exist far from the point $(1,0)$. 

In addition, we shall explain bifurcations of multiple kinks from 
$u_- = -\pi$ to $u_+ = \pi ( 2n - 1 )$, $n > 1$ in the discrete 
sine-Gordon equation from the point of normal form analysis. These 
bifurcations may occur along an infinite set of curves on the plane 
$(c,h)$ which all intersect the point $(1,0)$ (see numerical results 
in \cite{ACR03}). We shall also derive the truncated normal form 
from the discretizations of the inverse method reported in 
\cite{FZK99} and show that it may admit special solutions for 
monotonic kinks, when the continuity condition 
(\ref{nonlinearity-continuous}) is violated. 

We emphasize that our methods are very different from the
computations of beyond-all-order expansions, exploited in the 
context of difference maps in \cite{TTH98,T00} and differential 
advanced-delay equations in \cite{OPB05}. By working near the 
particular point $c = 1$ and $h = 0$, we avoid beyond-all-order 
expansions and derive non-existence results at the polynomial normal 
form. In an asymptotic limit of common validity, bifurcations of 
heteroclinic orbits in the truncated normal form can be studied with 
beyond-all-order computations (see recent analysis and review in 
\cite{TP05}). 

Our paper has the following structure. Section 2 discusses 
eigenvalues of the linearization problem at the zero equilibrium 
point and gives a formal derivation of the main result, the scalar 
fourth-order equation. The rigorous derivation of the scalar 
fourth-order equation from decompositions of solutions, projection 
techniques and truncation is described in Section 3. Section 4 
contains numerical analysis of the fourth-order equation, where we 
compute the split function for heteroclinic orbits. Persistence 
analysis of bounded solutions in the differential advance-delay 
equation is developed in Section 5. Section 6 concludes the paper 
with summary and open problems. A number of important but technical 
computations are reported in appendices to this paper. Appendix A 
contains the comparison of the fourth-order equation with the normal 
form from \cite{I95}. Appendix B gives the proof of existence of 
center manifold in the full system, which supplements the analysis 
in \cite{IK00}. Appendix C describes computations of the Stokes 
constant for heteroclinic orbits in the fourth-order differential 
equation from methods of \cite{TP05}. Appendix D discusses 
applications of the fourth-order equation to the inverse method from 
\cite{FZK99}.

\section{Resonances in dispersion relations and the scalar normal form}

When the differential advance-delay equation (\ref{advance-delay})
is linearized near the zero equilibrium point, we look for
solutions in the form $\phi(z) = e^{\lambda z}$, where $\lambda$
belongs to the set of eigenvalues. All eigenvalues of the
linearized problems are obtained from roots of the dispersion
relation:
\begin{equation}
\label{dispersion-relation} D(\Lambda; c,h) = 2 \left( \cosh
\Lambda - 1 \right) + h^2 - c^2 \Lambda^2 = 0,
\end{equation}
where $\Lambda = \lambda h$ is the eigenvalue in zoomed variable
$\zeta = \frac{z}{h}$. We are interested to know how many
eigenvalues occur on the imaginary axis and whether the imaginary
axis is isolated from the set of complex eigenvalues (that is
complex eigenvalues do not accumulate at the imaginary axis).
Imaginary eigenvalues $\Lambda = 2 i K$ of the dispersion relation
(\ref{dispersion-relation}) satisfy the transcendental equation:
\begin{equation}
\label{roots-dispersion} \sin^2 K = \frac{h^2}{4} + c^2 K^2.
\end{equation}
Particular results on roots of the transcendental equation
(\ref{roots-dispersion}) are easily deduced from analysis of the
function $\sin(x)$ in complex domain. Due to the obvious symmetry,
we shall only consider the non-negative values of $c$ and $h$.

\begin{itemize}
\item When $c = 0$ and $0 < h < 2$, equation
(\ref{roots-dispersion}) has only simple real roots $K$, such that
all eigenvalues $\Lambda$ are purely imaginary and simple. All
real roots $K$ become double at $h = 0$ and $h = 2$.

\item When $h = 0$ and $c \neq 1$, a double zero root of $K$ (a
double zero eigenvalue $\Lambda$) exists. When $0 < c < 1$, finitely
many purely imaginary eigenvalues $\Lambda$ exist (e.g. only one
pair of roots $K = \pm P$ exists for $c_0 < c < 1$, where $c_0 
\approx 0.22$), while infinitely many roots are complex and distant 
away the imaginary axis. When $c > 1$, all non-zero eigenvalues 
$\Lambda$ are complex and distant away the imaginary axis.
\end{itemize}

In the general case of $h \neq 0$ and $c \neq 0$, the
transcendental equation (\ref{roots-dispersion}) is more
complicated but it can be analyzed similarly to the dispersion
relation considered in \cite{I00,IK00}.

\begin{Lemma}
\label{lemma-eigenvalues} There exists a curve $h = h_*(c)$, which 
intersects the points $(1,0)$ and $(0,2)$ on the plane $(c,h)$, such 
that for $0 < c < 1$ and $0 < h < h_*(c)$ finitely many eigenvalues 
$\Lambda$ of the dispersion relation (\ref{dispersion-relation}) are 
located on the imaginary axis and all other eigenvalues are in a 
complex plane distant from the imaginary axis. The curve $h = 
h_*(c)$ corresponds to the 1:1 resonance Hopf bifurcation, where the 
set of imaginary eigenvalues includes only one pair of double 
eigenvalues.
\end{Lemma}

\begin{Remark}
The 1:1 resonant Hopf bifurcation is illustrated on Figure 1. Figure
1(a) shows a graphical solution of the transcendental equation 
(\ref{roots-dispersion}) for $h = h_*(c)$. Figure 1(b) shows the 
bifurcation curve $h = h_*(c)$ on the plane $(c,h)$. The curve $h = 
0$ for $0 < c < 1$ corresponds to the double zero eigenvalue in 
resonance with pairs of simple purely imaginary eigenvalues.
\end{Remark}

\begin{Proof}
Let $\Lambda = p + i q$ and rewrite the dispersion relation in the
equivalent form:
\begin{eqnarray*}
c^2 ( q^2 - p^2 ) + h^2 + 2(\cosh p \cos q - 1) & = & 0 \\
c^2 p q - \sinh p \sin q & = & 0.
\end{eqnarray*}
It follows from the system with $c \neq 0$ that the imaginary
parts of the eigenvalues is bounded by the real parts of the
eigenvalues:
\begin{equation}
q^2 \leq p^2 + \frac{4}{c^2} \cosh^2 \frac{p}{2}.
\end{equation}
Therefore, if complex eigenvalues accumulate to the imaginary axis, 
such that there exists a sequence $(p_n,q_n)$ with $\lim_{n \to 
\infty} p_n = 0$ and $\lim_{n \to \infty} q_n = q_*$, the 
accumulation point $(0,q_*)$ is bounded. However, since 
$D(\Lambda;c,h)$ is analytic in $\Lambda$, the dispersion relation 
(\ref{dispersion-relation}) may have finitely many roots of finite 
multiplicities in a bounded domain of the complex plane. Therefore, 
the accumulation point $(0,q_*)$ does not exist and complex 
eigenvalues are distant from the imaginary axis. By the same reason, 
there exist finitely many eigenvalues on the imaginary axis $\Lambda 
\in i \mathbb{R}$.

Looking at the double roots at $K = \pm P$, we find a
parametrization of the curve $h = h_*(c)$ in the form:
\begin{equation}
\label{reversible-bifurcation} c^2 = \frac{\sin P \cos P}{P},
\qquad h^2 = 4 \sin P \left( \sin P - P \cos P \right).
\end{equation}
A simple graphical analysis of the transcendental equation
(\ref{roots-dispersion}) shows that the double roots at $K = \pm
P$ are unique for $0 < P < \frac{\pi}{2}$, when the intersection
of the parabola and the trigonometric function occurs at the first
fundamental period of $\sin^2$-function. At $P=0$ (at the point
$(c,h) = (1,0)$), the pair of double imaginary eigenvalues
$\Lambda = \pm 2 i P$ merge at $\Lambda = 0$ and form a quadruple
zero eigenvalue. At $P = \frac{\pi}{2}$ (at the point $(c,h) =
(0,2)$), a sequence of infinitely many double imaginary
eigenvalues exists at $\Lambda = i \pi (1+2n)$, $n \in
\mathbb{Z}$.
\end{Proof}

\begin{figure}[htbp]
\begin{center}
\includegraphics[height=6cm]{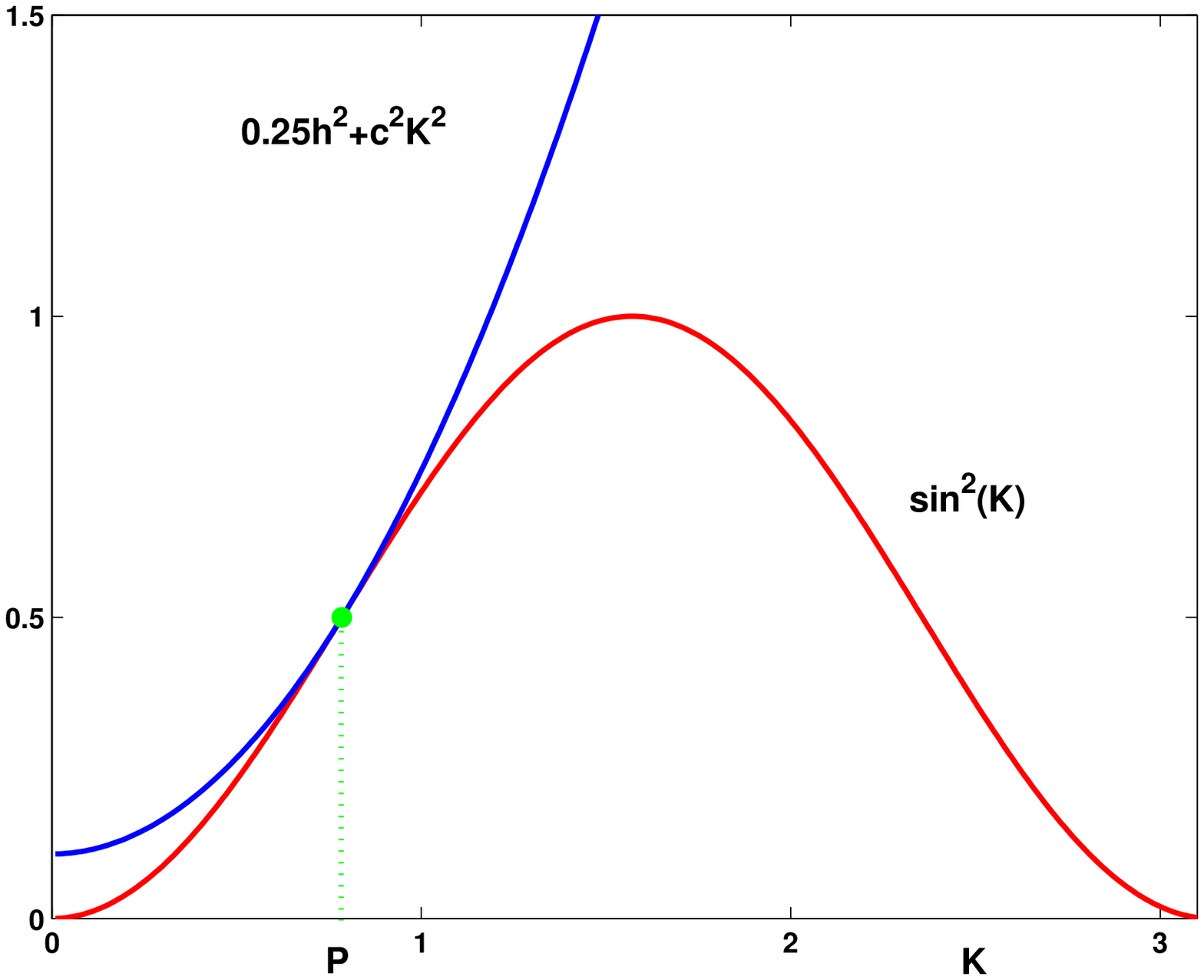}
\includegraphics[height=6cm]{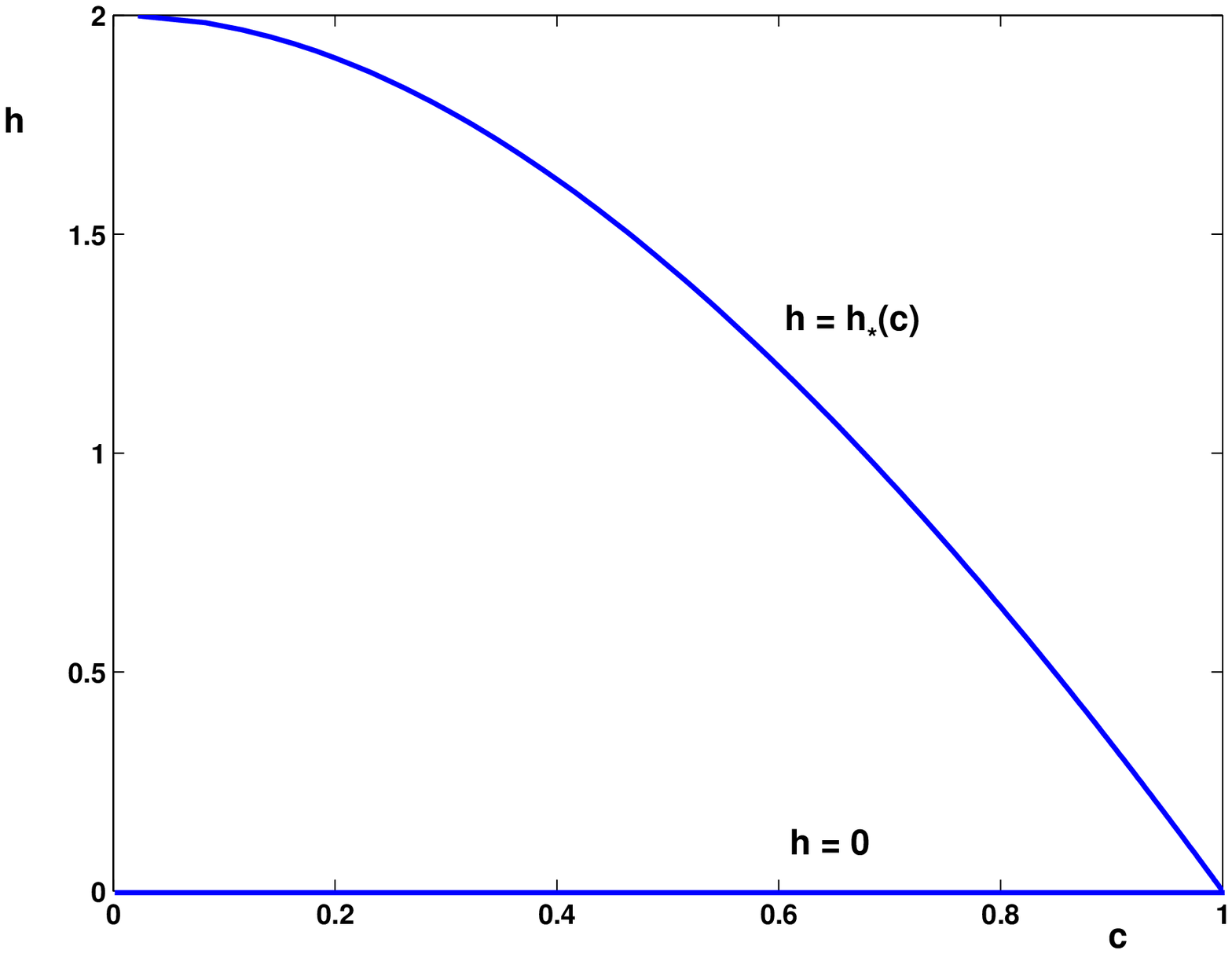}
\end{center}
\caption{Left: The graphical solution of the transcendental
equation (\ref{roots-dispersion}) for $P = \frac{\pi}{4}$, where
$(c,h)$ are given by (\ref{reversible-bifurcation}). Right: The
bifurcation curve for 1:1 resonant Hopf bifurcation $h = h_*(c)$
and double zero resonance bifurcation $h = 0$.}
\end{figure}

We will be interested in the reduction of the differential
advanced--delay equation (\ref{advance-delay}) at the particular
point $(c,h) = (1,0)$. Let $\epsilon$ be a small parameter that
defines a point on the plane $(c,h)$, which is locally close to
the bifurcation point $(1,0)$:
\begin{equation}
\label{perturbed-point} c^2 = 1 + \epsilon \gamma, \qquad h^2 =
\epsilon^2 \tau.
\end{equation}
Let $\Lambda = \sqrt{\epsilon} \Lambda_1$ be the scaled eigenvalue
near $\Lambda = 0$, so that the Taylor series expansion of the
dispersion relation (\ref{dispersion-relation}) allows for a
non-trivial balance of $\epsilon$-terms. Truncating the Taylor
series beyond the leading order ${\rm O}(\epsilon^2)$, we obtain
the bi-quadratic equation for the scaled eigenvalue $\Lambda_1$:
\begin{equation}
\label{quadric-curve} \frac{1}{12} \Lambda_1^4 - \gamma
\Lambda_1^2 + \tau = 0.
\end{equation}
The 1:1 resonance Hopf bifurcation corresponds to the point, where 
the roots of the bi-quadratic equation (\ref{quadric-curve}) are 
double and purely imaginary. This point occurs when $\tau = 3 
\gamma^2$, which agrees with the Taylor series expansion of the 
system (\ref{reversible-bifurcation}) in powers of $P$ at the 
leading orders of $P^2$ and $P^4$. As a result, the bifurcation 
curve $h = h_*(c)$ has the leading order behavior:
\begin{equation}
\label{h-star-c} h_*(c) = \sqrt{3} ( 1 - c^2) + {\rm o}(1-c^2).
\end{equation}
The normal form for travelling kinks, which is the main result of
our paper, can be recovered with a formal asymptotic expansion of
the nonlinear differential advance-delay equation
(\ref{advance-delay}). Let $\zeta = \frac{z}{h}$ and rewrite the
main equation (\ref{advance-delay}) in the form:
\begin{equation}
\label{advance-delay-Taylor} c^2 \phi''(\zeta) = \phi(\zeta + 1) -
2 \phi(\zeta) + \phi(\zeta - 1) + h^2 \phi(\zeta) + h^2
Q(\phi(\zeta - 1),\phi(\zeta),\phi(\zeta + 1)).
\end{equation}
Let $(c,h)$ be the perturbed bifurcation point in the form
(\ref{perturbed-point}). Let $\phi(\zeta)$ be a smooth function of
a slowly varying variable $\zeta_1 = \sqrt{\epsilon} \zeta$.
Expanding $\phi(\zeta_1 \pm \sqrt{\epsilon})$ in the Taylor series
in $\sqrt{\epsilon}$ and truncating by terms beyond the leading
order $\epsilon^2$, we obtain the scalar fourth-order differential
equation:
\begin{equation}
\label{scalar-fourth-order} \frac{1}{12} \phi^{({\rm
iv})}(\zeta_1) - \gamma \phi''(\zeta_1) + \tau F(\phi(\zeta_1)) =
0,
\end{equation}
where we have used the representations
(\ref{nonlinearity-continuous}) and 
(\ref{nonlinearity-representation}). The linearization of the 
nonlinear ODE (\ref{scalar-fourth-order}) near the zero equilibrium 
point recovers the bi-quadratic dispersion relation 
(\ref{quadric-curve}). The nonlinear equation 
(\ref{scalar-fourth-order}) has the equilibrium points $u_-$, $0$, 
and $u_+$, which are inherited from the equilibrium points of the 
Klein--Gordon equation (\ref{KG}). Our main examples will include 
the $\phi^4$ model with $F(u) = u (1 - u^2)$, $u_{\pm} = \pm 1$ and 
the sine--Gordon equation with $F(u) = \sin(u)$, $u_{\pm} = \pm 
\pi$.

\section{Decompositions, projections and truncation}

We shall derive the normal form equation (\ref{scalar-fourth-order}) 
with rigorous analysis when the solution of the differential 
advance-delay equation (\ref{advance-delay-Taylor}) is decomposed 
near a quadruple zero eigenvalue. We adopt notations of \cite{IK00} 
(see review in \cite{IJ04}) and rewrite the differential 
advance-delay equation as an infinite-dimensional evolution problem. 
We shall work with the scaled (inner) variable $\zeta = 
\frac{z}{h}$, where the differential advance-delay equation 
(\ref{advance-delay}) takes the form (\ref{advance-delay-Taylor}). 
Let $p$ be a new independent variable, such that $p \in [-1,1]$ and 
define the vector ${\bf U} = (U_1(\zeta),U_2(\zeta),U_3(\zeta,p) 
)^T$, such that
\begin{equation}
U_1 = \phi(\zeta), \qquad U_2 = \phi'(\zeta), \qquad U_3 =
\phi(\zeta + p).
\end{equation}
It is clear that $U_3(\zeta,p) = U_1(\zeta + p)$ and $U_3(\zeta,0)
= U_1(\zeta)$. The difference operators are then defined as
\begin{equation}
\label{difference-operators} U_3(\zeta,\pm 1) = \delta^{\pm} 
U_3(\zeta,p) = \phi(\zeta \pm 1) = \delta_{\pm} U_1(\zeta) = 
U_1(\zeta \pm 1).
\end{equation}
Let ${\cal D}$ and ${\cal H}$ be the following Banach spaces for 
${\bf U} = (U_1,U_2,U_3(p))^T$,
\begin{eqnarray}
\label{Banach1} {\cal D} & = & \left\{ (U_1,U_2) \in \mathbb{R}^2,
\; U_3 \in C^1([-1,1],\mathbb{R}),  \; U_3(0) = U_1 \right\}, \\
\label{Banach2} {\cal H} & = & \left\{ (U_1,U_2) \in \mathbb{R}^2,
\; U_3 \in C^0([-1,1],\mathbb{R}) \right\},
\end{eqnarray}
with the usual supremum norm. We look for a smooth mapping $\zeta
\mapsto {\bf U}(\zeta)$ in $C^0(\mathbb{R};{\cal D})$, which
represents classical solutions of the infinite-dimensional evolution 
problem:
\begin{equation}
\label{dynamical-system} \frac{d}{d\zeta} {\bf U} = {\cal L}_{c,h}
{\bf U} + \frac{h^2}{c^2} {\bf M}_0({\bf U}),
\end{equation}
where ${\cal L}_{c,h}$ and ${\bf M}_0({\bf U})$ are found from the
differential advance-delay equation (\ref{advance-delay-Taylor}):
\begin{equation}
\label{linear-operator} {\cal L}_{c,h} = \left( \begin{array}{ccc} 0 & 1 & 0 \\
\frac{h^2 - 2}{c^2} & 0 & \frac{1}{c^2} \left( \delta^+ + \delta^-
\right) \\ 0 & 0 & \partial_p \end{array} \right)
\end{equation}
and
\begin{equation}
\label{nonlinearity} {\bf M}_0({\bf U}) = (0,Q(\delta^-
U_3,U_1,\delta^+ U_3),0)^T.
\end{equation}
The linear operator ${\cal L}_{c,h}$ maps ${\cal D}$ into ${\cal H}$ 
continuously and it has a compact resolvent in ${\cal H}$. The 
nonlinearity ${\bf M}_0({\bf U})$ is analytic in an open 
neighborhood of ${\bf U} = {\bf 0}$, maps ${\cal H}$ into ${\cal D}$ 
continuously, and $\| {\bf M}_0({\bf U}) \|_{{\cal D}} = {\rm 
O}\left(\| {\bf U} \|_{{\cal H}}^2\right)$.

The spectrum of ${\cal L}_{c,h}$ consists of an infinite set of
isolated eigenvalues of finite multiplicities. By virtue of the
Laplace transform, eigenvalues of the linear operator ${\cal
L}_{c,h}$ are found with the solution ${\bf U}(\zeta,p) = (1,
\Lambda, e^{\Lambda p})^T e^{\Lambda \zeta}$, when the linear
problem ${\bf U}'(\zeta) = {\cal L}_{c,h} {\bf U}$ is reduced to the 
dispersion relation (\ref{dispersion-relation}), that is 
$D(\Lambda;c,h) = 0$. We are particularly interested in the 
bifurcation point $c = 1$ and $h = 0$, when
\begin{equation}
\label{dispersion-relation-bifurcation} D_1( \Lambda) \equiv
D(\Lambda;1,0) = 2 (\cosh \Lambda - 1) - \Lambda^2.
\end{equation}
The transcendental equation $D_1(\Lambda) = 0$ has the quadruple 
zero root and no other root in the neighborhood of $\Lambda \in i 
\mathbb{R}$. The four generalized eigenvectors of the Jordan chain 
for the zero eigenvalue, ${\cal L}_{1,0} {\bf U}_j = {\bf U}_{j-1}$, 
$j=1,2,3,4$ with ${\bf U}_{-1} = {\bf 0}$, are found exactly as:
\begin{equation}
\label{zero-eigenvectors} {\bf U}_0 = \left( \begin{array}{cc} 1
\\ 0 \\ 1 \end{array} \right), \qquad {\bf U}_1 = \left( \begin{array}{cc}
0 \\ 1 \\ p \end{array} \right), \qquad {\bf U}_2 = \left(
\begin{array}{cc} 0 \\ 0 \\ \frac{1}{2} p^2 \end{array} \right), \qquad {\bf U}_3
= \left( \begin{array}{cc} 0 \\ 0 \\ \frac{1}{6} p^3 \end{array}
\right).
\end{equation}

The dynamical system (\ref{dynamical-system}) has the
reversibility symmetry ${\cal S}$, such that
\begin{equation}
\label{dynamical-system-reversibility} -\frac{d}{d\zeta} {\cal S}
{\bf U} = {\cal L}_{c,h} {\cal S} {\bf U} + \frac{h^2}{c^2} {\bf
M}_0({\cal S} {\bf U}),
\end{equation}
where
\begin{equation}
\label{reversibility} {\cal S} {\bf U} = (U_1(\zeta), -U_2(\zeta),
U_3(\zeta,-p) )^T
\end{equation}
and the symmetry property (\ref{nonlinearity-symmetry}) has been 
used. Since both the linear operator (\ref{linear-operator}) and the 
nonlinearity (\ref{nonlinearity}) anti-commute with the 
reversibility operator (\ref{reversibility}), the standard property 
of reversible systems holds: if ${\bf U}(\zeta)$ is a solution of 
(\ref{dynamical-system}) for forward time $\zeta > 0$, then ${\cal 
S} {\bf U}(-\zeta)$ is a solution of (\ref{dynamical-system}) for 
backward time $\zeta < 0$ (see \cite{LR98} for review). Applying the 
reversibility operator ${\cal S}$ to the eigenvectors 
(\ref{zero-eigenvectors}) at the bifurcation of the quadruple zero 
eigenvalue, we observe that
\begin{equation}
\label{reversibility-eigenvectors} {\cal S} {\bf U}_j = (-1)^j
{\bf U}_j, \quad j = 0,1,2,3.
\end{equation}
This bifurcation case fits to the analysis of \cite{I95}, where the 
reversible normal form was derived for the quadruple zero 
eigenvalue. Appendix A shows that the general reversible normal form 
from \cite{I95} is reduced to the normal form equation 
(\ref{scalar-fourth-order}) by appropriate scaling. However, we 
notice that this result cannot be applied directly, since it is only 
valid in a neighborhood of the origin, while in the present case, 
the solutions we are interested in are of order of ${\rm O}(1)$. 

In order to study bifurcation of quadruple zero eigenvalue in the 
reversible system (\ref{dynamical-system}), we shall define the 
perturbed point $(c,h)$ near the bifurcation point $(1,0)$ with an 
explicit small parameter $\varepsilon$. Contrary to the definition 
(\ref{perturbed-point}), it will be easier to work with the 
parameters $(\gamma,\tau)$, defined from the relations:
\begin{equation}
\label{perturbed-point-bifurcation} \frac{1}{c^2} = 1 -
\varepsilon \gamma, \qquad \frac{h^2}{c^2} = \varepsilon^2 \tau,
\end{equation}
where $\varepsilon$ is different from $\epsilon$ used in
(\ref{perturbed-point}). The dynamical system
(\ref{dynamical-system}) with the parametrization
(\ref{perturbed-point-bifurcation}) is rewritten in the explicit
form:
\begin{equation}
\label{dynamical-system-bifurcation} \frac{d}{d\zeta} {\bf U} =
{\cal L}_{1,0} {\bf U} + \varepsilon \gamma {\bf L}_1({\bf U}) +
\varepsilon^2 \tau {\bf L}_2({\bf U}) + \varepsilon^2 \tau {\bf
M}_0({\bf U}),
\end{equation}
where ${\cal L}_{1,0}$ follows from (\ref{linear-operator}) with
$(c,h) = (1,0)$, ${\bf M}_0({\bf U})$ is the same as
(\ref{nonlinearity}), and the linear terms ${\bf L}_{1,2}({\bf
U})$ are
\begin{equation}
\label{perturbations} {\bf L}_1({\bf U}) = (0,2 U_1 - (\delta^+ +
\delta^-) U_3,0)^T, \qquad {\bf L}_2({\bf U}) = (0,U_1,0)^T.
\end{equation}
Bounded solutions of the ill-posed initial-value problem 
(\ref{dynamical-system-bifurcation}) on the entire axis $\zeta \in 
\mathbb{R}$ are subject of our interest, with the particular 
emphasis on the kink solutions (heteroclinic orbits between non-zero 
equilibrium points $u_-$ and $u_+$). These bounded solutions can be 
constructed with the decomposition of the solution of the 
infinite-dimensional system (\ref{dynamical-system-bifurcation}), 
projection to the finite-dimensional subspace of zero eigenvalue and 
an infinite-dimensional subspace of non-zero eigenvalues, and 
truncation of the resulting system of equations into the 
fourth-order differential equation (\ref{scalar-fourth-order}). 
Since the original system is reversible, the fourth-order equation 
must inherit the reversibility property.

The techniques of decompositions, projections and truncation rely
on the solution of the resolvent equation $(\Lambda {\cal I} -
{\cal L}_{1,0}) {\bf U} = {\bf F}$, where ${\bf U} \in {\cal D}$,
${\bf F} \in {\cal H}$ and $\Lambda \in \mathbb{C}$. When
$\Lambda$ is different from the roots of $D_1(\Lambda)$, the
explicit solution of the resolvent equation is obtained in the
form:
\begin{eqnarray}
\nonumber
U_1 & = & - \left[ D_1(\Lambda) \right]^{-1} \tilde{F}(\Lambda), \\
\nonumber
U_2 & = & \Lambda U_1 - F_1, \\
\label{solution-resolvent-equation-0} U_3(p) & = & U_1 e^{\Lambda p} 
- \int_0^p F_3(s) e^{\Lambda (p-s)} ds,
\end{eqnarray}
where
\begin{equation}
\label{solution-resolvent-equation} \tilde{F}(\Lambda) = \Lambda F_1 
+ F_2 - \int_0^1 \left( F_3(s) e^{\Lambda (1-s)} - F_3(-s) 
e^{-\Lambda (1-s)} \right) ds.
\end{equation}
The quadruple zero eigenvalue appears as the quadruple pole in the
solution ${\bf U}(p;\Lambda)$ of the resolvent equation near
$\Lambda = 0$. Let us decompose the solution of the dynamical system 
(\ref{dynamical-system-bifurcation}) into two parts:
\begin{equation}
\label{decomposition1} {\bf U} = {\bf X} + {\bf Y},
\end{equation}
where ${\bf X}$ is a projection to the fourth-dimensional subspace 
of the quadruple zero eigenvalue,
\begin{equation}
\label{decomposition2} {\bf X} = A(\zeta) {\bf U}_0 + B(\zeta) {\bf 
U}_1(p) + C(\zeta) {\bf U}_2(p) + D(\zeta) {\bf U}_3(p)
\end{equation}
and ${\bf Y}$ is the projection on the complementary invariant 
subspace of ${\cal L}_{1,0}$. We notice in particular that 
\begin{eqnarray}
\label{U1}
U_{1}(\zeta ) &=& A(\zeta )+Y_{1}(\zeta ) \\
\label{U3}
U_{3}(\zeta ,p) &=&A(\zeta )+pB(\zeta )+\frac{1}{2}p^{2}C(\zeta )+\frac{1}{6}%
p^{3}D(\zeta )+Y_{3}(\zeta,p).
\end{eqnarray}
With the decomposition 
(\ref{decomposition1})--(\ref{decomposition2}), the problem 
(\ref{dynamical-system-bifurcation}) is rewritten in the form:
\begin{equation}
\label{resolvent-equation} \frac{d}{d\zeta} {\bf Y} - {\cal
L}_{1,0} {\bf Y} = {\bf F}(A,B,C,D,{\bf Y}),
\end{equation}
where
\begin{eqnarray*}
F_1 & = & B - A' \\
F_2 & = & C - B' - \varepsilon \gamma \left( C - 2 Y_1 + (\delta^+ + 
\delta^-) Y_3 \right) + \varepsilon^2 \tau \left( A + Q + Y_1 \right) \\
F_3 & = & B - A' + p(C - B') + \frac{1}{2} p^2 (D - C') +
\frac{1}{6} p^3 (-D'),
\end{eqnarray*}
and
$$
Q = Q\left( A - B + \frac{1}{2} C -\frac{1}{6} D + \delta^- Y_3, A
+ Y_1, A + B + \frac{1}{2} C +\frac{1}{6} D + \delta^+ Y_3\right).
$$
We notice that $p$ is a dumb variable in (\ref{resolvent-equation}). 
Since the vector ${\bf Y}$ is projected to the subspace of non-zero 
eigenvalues of ${\cal L}_{1,0}$, the right-hand-side of the 
"inhomogeneous" problem (\ref{resolvent-equation}) must be 
orthogonal to eigenvectors of the adjoint operator ${\cal 
L}_{1,0}^*$ associated to the zero eigenvalue. Equivalently, one can 
consider the solution of the resolvent equation 
(\ref{solution-resolvent-equation}) and remove the pole 
singularities from the function $\left[ D_1(\Lambda) \right]^{-1} 
\tilde{F}(\Lambda)$ near $\Lambda = 0$. See p.853 in \cite{I00} for 
the projections on the four-dimensional subspace and, in particular, 
the formula:
\begin{equation}
\label{F-0} {\cal P}_h \left( \begin{array}{cc} 0 \\ -1 \\ 0 
\end{array} \right) = \left( \begin{array}{cc} 0 \\ -\frac{3}{5} \\ 
\frac{2}{5} p(1 - 5 p^2) \end{array} \right) \equiv {\bf F}_0(p),
\end{equation}
where ${\cal P}_h$ is a projection operator to the space for ${\bf 
Y}$. This procedure results in the fourth-order differential system 
for components $(A,B,C,D) \in \mathbb{R}^4$ in $\zeta$:
\begin{eqnarray}
\nonumber A' & = & B \\
\nonumber B' & = & C - \frac{2}{5} \varepsilon \gamma \left( C - 2
Y_1 + (\delta^+ + \delta^-) Y_3 \right) + \frac{2}{5}
\varepsilon^2 \tau \left( A + Q + Y_1 \right) \\
\nonumber C' & = & D \\
\label{center-manifold} D' & = & 12 \varepsilon \gamma \left( C -
2 Y_1 + (\delta^+ + \delta^-) Y_3 \right) - 12 \varepsilon^2 \tau
\left( A + Q + Y_1 \right).
\end{eqnarray}
While the system (\ref{resolvent-equation}) has the reversibility 
symmetry ${\cal S}$ in (\ref{reversibility}), the system  
(\ref{center-manifold}) has the reduced reversibility symmetry,
\begin{equation}
\label{reversibility-reduced} {\cal S}_{0} {\bf A} = 
(A,-B,C,-D)^{T},
\end{equation}
where ${\bf A} = (A,B,C,D)^T$. The closed system of equations 
(\ref{resolvent-equation}) and (\ref{center-manifold}) can be 
rewritten by using the scaling transformation:
\begin{equation}
\label{scaling-transformation} A = \phi_1(\zeta_1), \quad B =
\sqrt{\varepsilon} \phi_2(\zeta_1), \quad C = \varepsilon
\phi_3(\zeta_1), \quad D = \varepsilon^{3/2} \phi_4(\zeta_1), \quad 
{\bf Y} = \varepsilon^{3/2} \mbox{\boldmath $\psi$}(\zeta,p),
\end{equation}
where $\zeta_1 = \sqrt{\varepsilon} \zeta$. The four-dimensional
vector $\mbox{\boldmath $\phi$} \in \mathbb{R}^4$ satisfies the
vector system:
\begin{equation}  
\label{center-manifold-Final} \frac{d}{d \zeta_1} \left( 
\begin{array}{cc} \phi_1 \\ \phi_2 \\ \phi_3 \\ \phi_4 \end{array} 
\right) = \left( \begin{array}{ccccc} 0 & 1 & 0 & 0  \\ 0 & 0 & 1 & 
0 \\ 0 & 0 & 0 & 1 \\ 0 & 0 & 0 & 0 
\end{array} \right) \left( 
\begin{array}{cc} \phi_1 \\ \phi_2 \\ \phi_3 \\ \phi_4 \end{array} 
\right) + \left( \begin{array}{cc} 0 \\ - \frac{2}{5} \varepsilon \\ 
0  \\ 12 \end{array} \right) g,
\end{equation}
where $g$ is a scalar function given by
$$
g = \gamma \left( \phi_3 - 2 \sqrt{\varepsilon} \psi_1 + 
\sqrt{\varepsilon} (\delta^+ + \delta^-) \psi_3 \right) - \tau 
\left( \phi_1 + Q + \varepsilon^{3/2} \psi_1 \right)
$$
and     
\begin{eqnarray*}
Q = Q\left( \phi_1 - \sqrt{\varepsilon} \phi_2 + \frac{1}{2} 
\varepsilon \phi_3 - \frac{1}{6}\varepsilon^{3/2} \phi_4 + 
\sqrt{\varepsilon^3} \delta^- \psi_3, \phi_1 + \varepsilon^{3/2} 
\psi_1, \right. & \\
\left. \phi_1 + \sqrt{\varepsilon} \phi_2 + \frac{1}{2} \varepsilon 
\phi_3 + \frac{1}{6} \sqrt{\varepsilon^3} \phi_4 + \varepsilon^{3/2} 
\delta^+ \psi_3 \right). & 
\end{eqnarray*}
The infinite-dimensional vector $\mbox{\boldmath $\psi$} \in {\cal 
D}$ satisfies the "inhomogeneous" problem:
\begin{equation}
\label{resolvent-equation-Final} \frac{d}{d\zeta} \mbox{\boldmath 
$\psi$} - {\cal L}_{1,0} \mbox{\boldmath $\psi$} = \sqrt{\epsilon} g 
{\bf F}_0(p),
\end{equation}
where ${\bf F}_0$ is given in (\ref{F-0}). The truncation of the 
fourth-order system (\ref{center-manifold-Final}) at $\varepsilon = 
0$ recovers the scalar fourth-order equation 
(\ref{scalar-fourth-order}). The set of equilibrium points of the 
coupled system 
(\ref{center-manifold-Final})--(\ref{resolvent-equation-Final}) 
contains the set of equilibrium points of the fourth-order equation 
(\ref{scalar-fourth-order}), since the constraints $\phi_2 = \phi_3 
= \phi_4 = 0$ and $\mbox{\boldmath $\psi$} = {\bf 0}$ reduce the 
coupled system to the algebraic equation $\phi_1 + 
Q(\phi_1,\phi_1,\phi_1)$, that is $F(\phi_1) = 0$. Our central 
result is the proof of existence of the center manifold in the 
coupled system 
(\ref{center-manifold-Final})--(\ref{resolvent-equation-Final}).

\begin{Theorem} 
Fix $M>0$ and let $\varepsilon >0$ be small enough. Then for any 
given $\mbox{\boldmath $\phi$} \in 
C_{b}^{0}(\mathbb{R},\mathbb{R}^{4})$, such that
\begin{equation}
||\mbox{\boldmath $\phi$}|| \leq M,
\end{equation}
there exists a unique solution $\mbox{\boldmath $\psi$} \in 
C_{b}^{0}(\mathbb{R},{\cal D})$ of the system 
(\ref{resolvent-equation-Final}), such that 
\begin{equation}
\| \mbox{\boldmath $\psi$} \| \leq \sqrt{\varepsilon} K,
\end{equation}
where the constant $K > 0$ is independent of $\varepsilon$. 
Moreover, when $\mbox{\boldmath $\phi$}(\zeta_1)$ is anti-reversible 
with respect to ${\cal S}_0$, i.e. ${\cal S}_0 \mbox{\boldmath 
$\phi$}(-\zeta_{1}) = -\mbox{\boldmath $\phi$}(\zeta_{1})$, then  
$\mbox{\boldmath $\psi$}(\zeta)$ is anti-reversible with respect to 
${\cal S}$, i.e. ${\cal S} \mbox{\boldmath $\psi$}(-\zeta) = 
-\mbox{\boldmath $\psi$}(\zeta)$. If the solution $\mbox{\boldmath 
$\phi$}(\zeta_1)$ tends towards a $T$-periodic orbit as $\zeta 
_{1}\rightarrow \pm \infty$, then the solution $\mbox{\boldmath 
$\psi$}(\zeta)$ tends towards a $T/\sqrt{\varepsilon}$-periodic 
orbit of the system (\ref{resolvent-equation-Final}) as $\zeta 
\rightarrow \pm \infty$. \label{theorem-center-manifold}
\end{Theorem}

In this theorem, $C_b^0(\mathbb{R},{\cal E})$ denotes the Banach 
space (sup norm) of continuous and bounded functions on $\mathbb{R}$ 
taking values in the Banach space ${\cal E}$. 

The proof of existence of the center manifold for the class of 
reversible dynamical system (\ref{dynamical-system}) is developed in 
Lemmas 2--4 of \cite{IK00} with the Green function technique. 
Applications of the center manifold reduction in a similar content 
can be found in \cite{I00,IJ04,FR05}. We note that the case of the 
quadruple zero eigenvalue was not studied in any of the previous 
publications (e.g. the pioneer paper \cite{IK00} dealt with the case 
of two pairs of double and simple imaginary eigenvalues). We develop 
an alternative proof of existence of center manifold in the coupled 
system 
(\ref{center-manifold-Final})--(\ref{resolvent-equation-Final}) with 
a more explicit and elementary method in Appendix B.

\section{Numerical analysis of the truncated normal form}

Eigenvalues of the linear part of the fourth-order equation 
(\ref{scalar-fourth-order}), that are roots of the bi-quadratic 
equation (\ref{quadric-curve}), are shown on Figure 2. Due to the 
definition (\ref{perturbed-point-bifurcation}), it makes sense to 
consider only the upper half plane of the parameter plane 
$(\gamma,\tau)$. The curve $\Gamma_1 = \{ \gamma < 0, \tau = 0 \}$ 
corresponds to the double zero eigenvalue that exists in resonance 
with a pair of purely imaginary eigenvalues. The curve $\Gamma_2 = 
\{ \gamma < 0, \tau = 3 \gamma^2 \}$ corresponds to the reversible 
1:1 resonance Hopf bifurcation $h = h_*(c)$ in the limit 
(\ref{h-star-c}). These two curves confine the parameter domain of 
our interest, where a heteroclinic orbit between two non-zero 
equilibrium points $u_-$ and $u_+$ of the nonlinearity function 
$F(u)$ may undertake a bifurcation. 

\begin{figure}[htbp]
\begin{center}
\includegraphics[height=8cm]{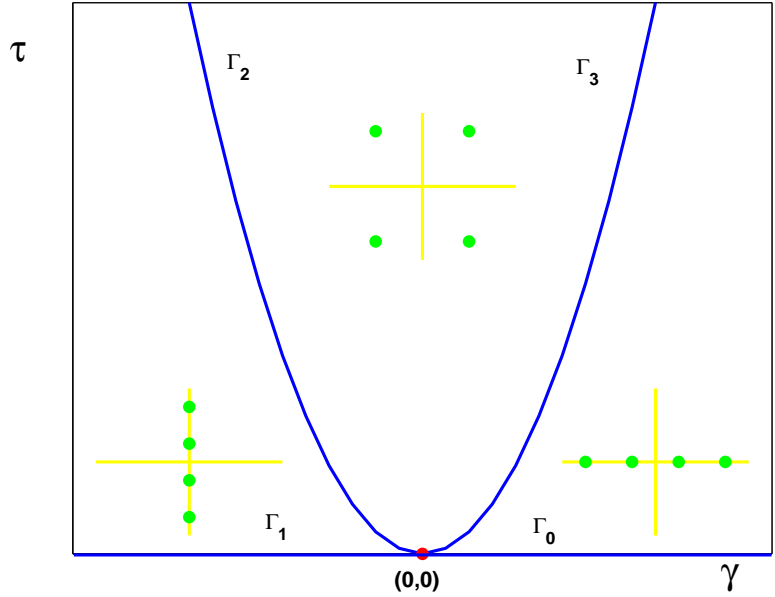}
\end{center}
\caption{Eigenvalues of the scalar normal form
(\ref{scalar-fourth-order}) on the parameter plane $(\gamma,\tau)$. 
The curves $\Gamma_{0,1,2,3}$ mark various bifurcations of 
co-dimension one. }
\end{figure}

A complete classification of various bifurcations of co-dimension 
one in the normal form (\ref{scalar-fourth-order}) can be found in 
\cite{I95}. Homoclinic orbits may bifurcate above the curve 
$\Gamma_2$ and may exist in the domain $D = \{ \tau 
> 3 \gamma^2, \; \gamma < 0 \} \cup \{ \tau > 0, \; \gamma \geq 0 
\}$. Homoclinic orbits correspond to localized solutions to the zero 
equilibrium state and they are not studied in this paper. We focus 
hence on bounded heteroclinic orbits which connect the non-zero 
equilibrium states $u_-$ and $u_+$.

Let us rewrite the scalar fourth-order equation
(\ref{scalar-fourth-order}) in the normalized form,
\begin{equation}
\label{ODE} \phi^{({\rm iv})} + \sigma \phi'' + F(\phi) = 0,
\end{equation}
where $\phi = \phi(t)$, $t = (12 \tau)^{1/4} \zeta_1$, and $\sigma
= -\frac{\sqrt{12} \gamma}{\sqrt{\tau}}$. The domain between the
curves $\Gamma_2$ and $\Gamma_1$ corresponds to the semi-infinite
interval $\sigma > 2$. We shall consider separately the discrete
$\phi^4$ model with $F(\phi) = \phi ( 1- \phi^2)$ and the
sine--Gordon model with $F(\phi) = \sin(\phi)$.

When $F(\phi) = \phi ( 1- \phi^2)$, there are two non-zero
equilibrium points $u_{\pm} = \pm 1$. Linearization of the scalar
equation (\ref{ODE}) near non-zero equilibrium points gives a pair
of real eigenvalues $(\lambda_0,-\lambda_0)$ and a pair of purely
imaginary eigenvalues $(i \omega_0,-i\omega_0)$ for any $\sigma
\in \mathbb{R}$. Our numerical algorithm is based on the
consideration of the unstable solution $\phi_u(t)$, such that
\begin{equation}
\lim_{t \to -\infty} \phi_u(t) = -1, \qquad \lim_{t \to -\infty}
\left( \phi_u(t) + 1 \right) e^{-\lambda_0 t} = c_0,
\end{equation}
where $c_0$ is arbitrary positive constant, which is related to the 
translation of the unstable solution. Let $t_0(c_0)$ be a zero of 
the unstable solution $\phi_u(t)$, if it exists. Since the ODE 
problem (\ref{ODE}) with odd function $F(-\phi) = -F(\phi)$ has the 
reflection symmetry $\phi \mapsto -\phi$, we look for a bounded 
heteroclinic orbit as an odd function of $t$. Therefore, the 
distance between stable and unstable solutions can be measured by 
the split function $K = \phi_u''(t_0(c_0))$, such that the roots of 
$K$ give odd heteroclinic orbits. Although the split function $K$ is 
computed at $t_0$, which depends on $c_0$, it follows from the 
translational invariance of the problem (\ref{ODE}) that $K$ is 
independent of $c_0$. (It may depend on the external parameter 
$\sigma$ of the problem (\ref{ODE}).) The first zero $t_0$ 
corresponds to the monotonic heteroclinic orbit when $K = 0$, while 
subsequent zeros correspond to non-monotonic heteroclinic orbits. We 
shall here consider monotonic heteroclinic orbits.

The unstable solution of the scalar equation (\ref{ODE}) is
numerically approximated by the solution of the initial-value 
problem with the initial data:
\begin{equation}
\label{initial-value} \left( \begin{array}{cc} \phi(0) \\ \phi'(0) \\
\phi''(0) \\ \phi'''(0) \end{array} \right) = \left( \begin{array}{cc} -1 \\ 0 \\
0 \\ 0 \end{array} \right) + c_0 \left( \begin{array}{cc} 1 \\
\lambda_0 \\\lambda_0^2 \\ \lambda_0^3 \end{array} \right),
\end{equation}
where $\lambda_0 > 0$ and $0 < c_0 \ll 1$. The solution of the 
initial-value problem is stopped at the first value $t_0 > 0$, where 
$\phi(t_0) = 0$ and the split function $K$ is approximated as $K = 
\phi''(t_0)$. Besides the external parameter $\sigma$, the numerical 
approximation of the split function $K$ may depend on the shooting 
parameter $c_0$ and the discretization parameter $dt$ of the ODE 
solver, such that $K = K(\sigma;c_0,dt)$.

The graph of $K$ versus $\sigma$ for a fixed set of values of $c_0$ 
and $dt$ is shown on Figure 3(a). It is clearly seen that the split 
function $K(\sigma)$ is non-zero for finite values of $\sigma$. It 
becomes exponentially small in $\sigma^{-1}$ as $\sigma \to \infty$. 
In the beyond-all-order asymptotic limit, computations of Stokes 
constants give equivalent information to the computation of the 
split function $K(\sigma)$ (see \cite{TTH98,T00}). We apply this 
technique in Appendix C to show that the Stokes constant for the 
scalar fourth-order equation (\ref{ODE})  is non-zero in the limit 
$\sigma \to \infty$.

Behavior of the split function $K$ for a fixed value of $\sigma$ is 
analyzed on Figures 4(a,b) as the values of $dt$ and $c_0$ are 
reduced. When $dt$ is sufficiently reduced, the value of $K$ 
converges to a constant value monotonically (see Figure 4(a)). When 
$c_0$ changes, the value of $K$ oscillates near a constant level. 
The amplitude of oscillations depend on the value of $dt$ and it 
becomes smaller when the value of $dt$ is reduced (see Figure 4(b)). 
Oscillations on Figure 4(b) can be used as the tool to measure the 
error of numerical approximations for the split function 
$K(\sigma)$. The maximum relative error is computed from the ratio 
of the standard deviation of the values of $K(\sigma)$ to the mean 
value when the parameter $c_0$ varies. The graph of the relative 
error is plotted on Figure 3(b) versus parameter $\sigma$. The 
relative error increases for larger values of $\sigma$ since 
$K(\sigma)$ becomes exponentially small in the limit $\sigma \to 
\infty$. 

\begin{figure}[htbp]
\begin{center}
\includegraphics[height=5.5cm]{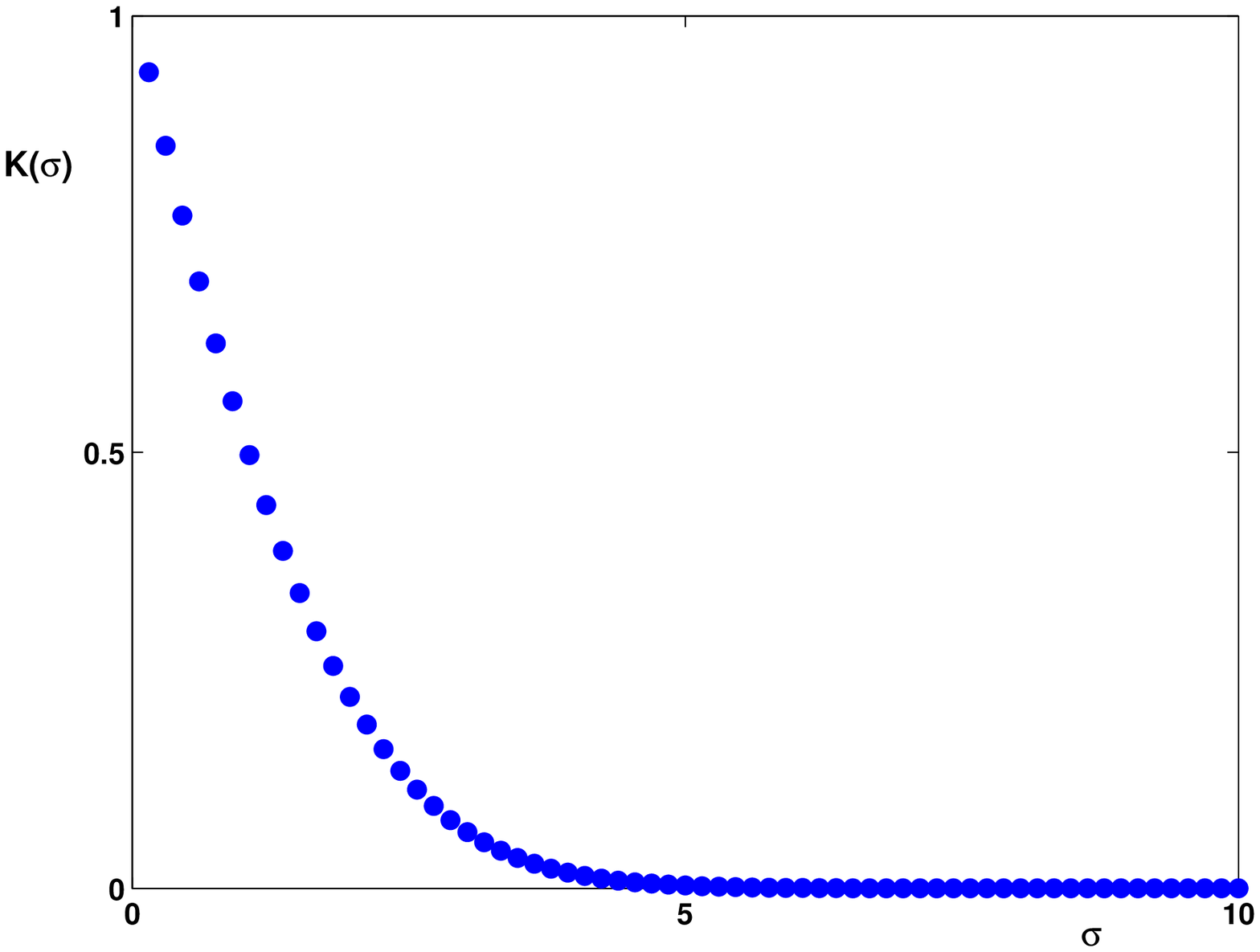} 
\includegraphics[height=5.5cm]{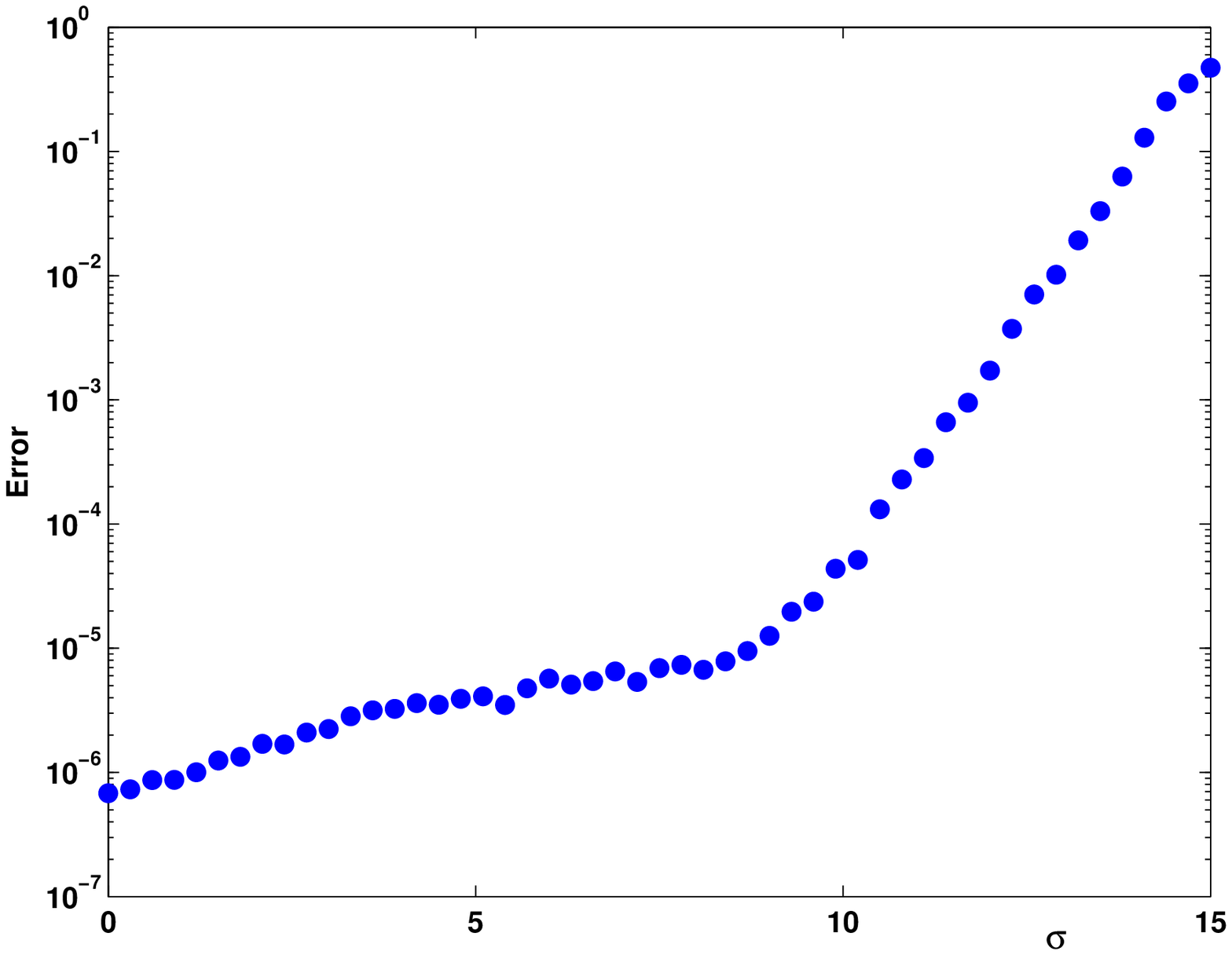}  
\end{center}
\caption{Left: The graph of the split function $K(\sigma)$ for $c_0 
= 0.00001$ and $dt = 0.005$ for kink solutions of the ODE 
(\ref{ODE}) with $F(\phi) = \phi (1 - \phi^2)$. Right: relative 
error of numerical approximations of $K(\sigma)$ for $dt = 0.005$ 
under variations of $c_0$.}
\end{figure}

\begin{figure}[htbp]
\begin{center}
\includegraphics[height=6cm]{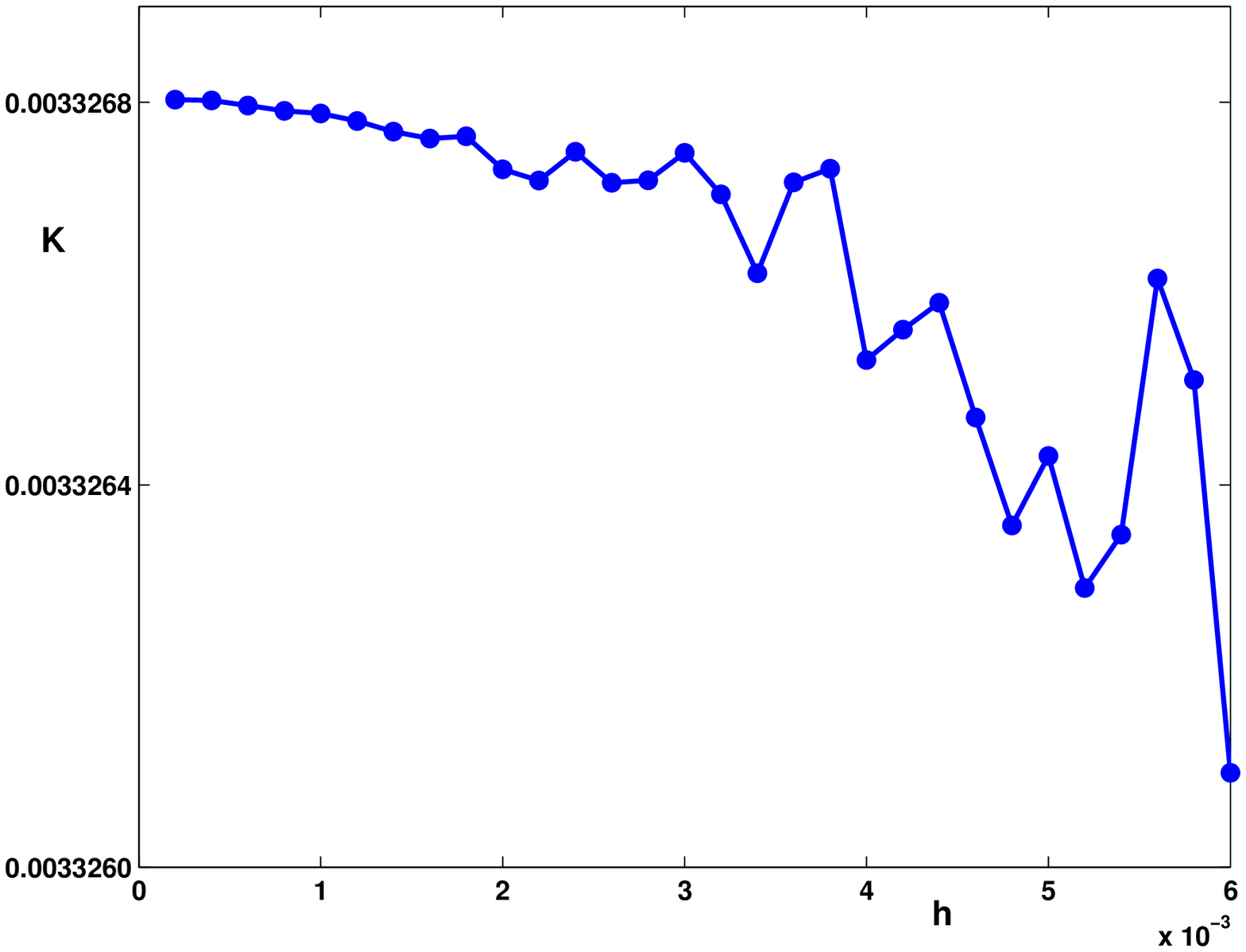}
\includegraphics[height=6cm]{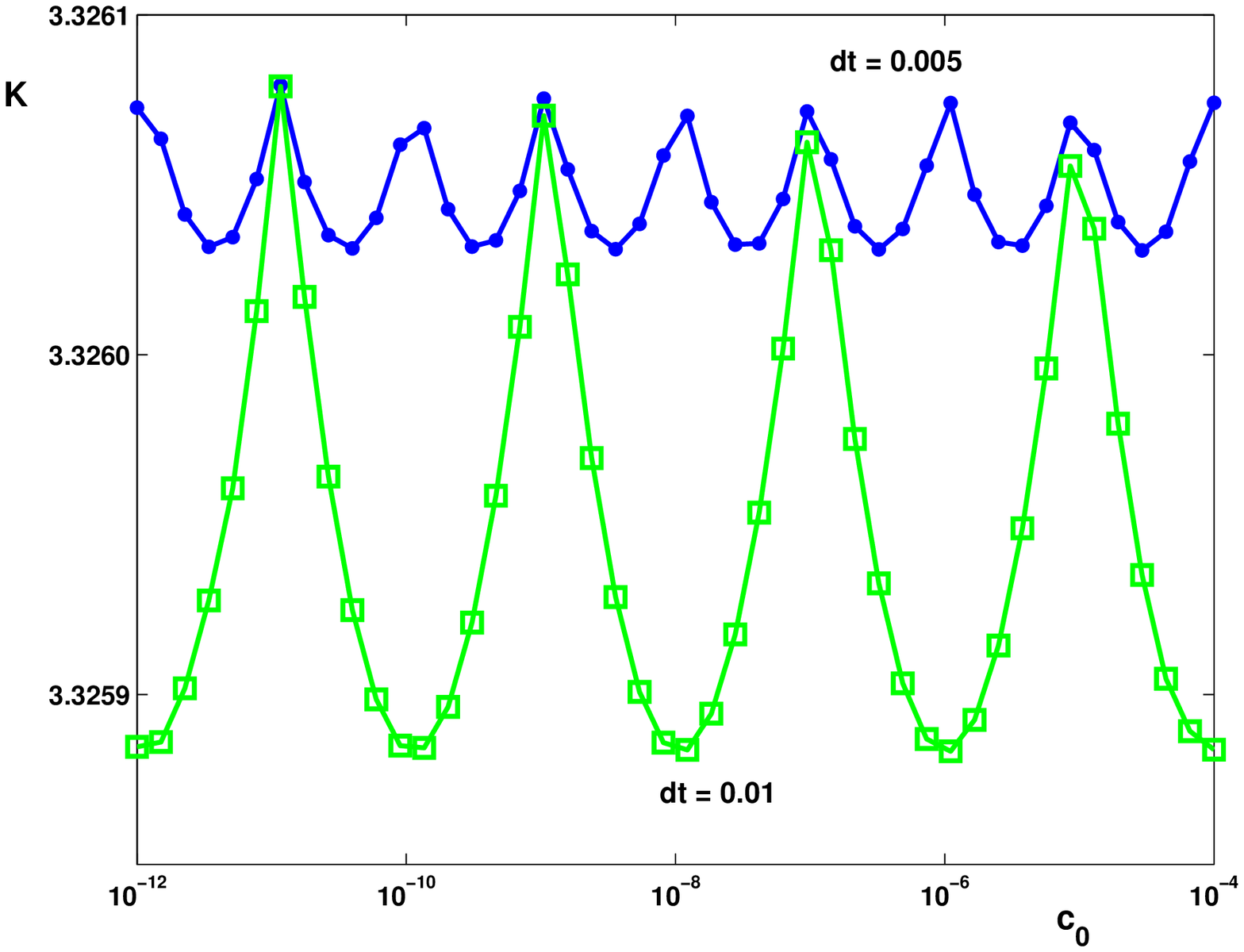}
\end{center}
\caption{Left: Convergence of the split function $K(\sigma)$ as $dt 
\to 0$ for $\sigma = 5$ and $c_0 = 0.00001$. Right: Oscillations of 
the split function $K(\sigma)$ as $c_0 \to 0$ for $\sigma = 5$ and 
$dt = 0.005$. }
\end{figure}

When $F(\phi) = \sin(\phi)$, there are two non-zero equilibrium
points $u_{\pm} = \pm \pi$. Linearization of the scalar equation
(\ref{ODE}) near non-zero equilibrium points still has a pair of
real eigenvalues $(\lambda_0,-\lambda_0)$ and a pair of purely
imaginary eigenvalues $(i \omega_0,-i\omega_0)$ for any $\sigma \in 
\mathbb{R}$. Therefore, we adjust the same algorithm to the 
sine--Gordon model. Figures 5(a) displays the graph of the split 
function $K$ versus $\sigma$ for a fixed set of values of $c_0$ and 
$dt$, which Figure 5(b) shows the maximum relative error of 
numerical approximations. The split function $K(\sigma)$ displays 
the same behavior as that for the $\phi^4$ model.

\begin{figure}[htbp]
\begin{center}
\includegraphics[height=5.5cm]{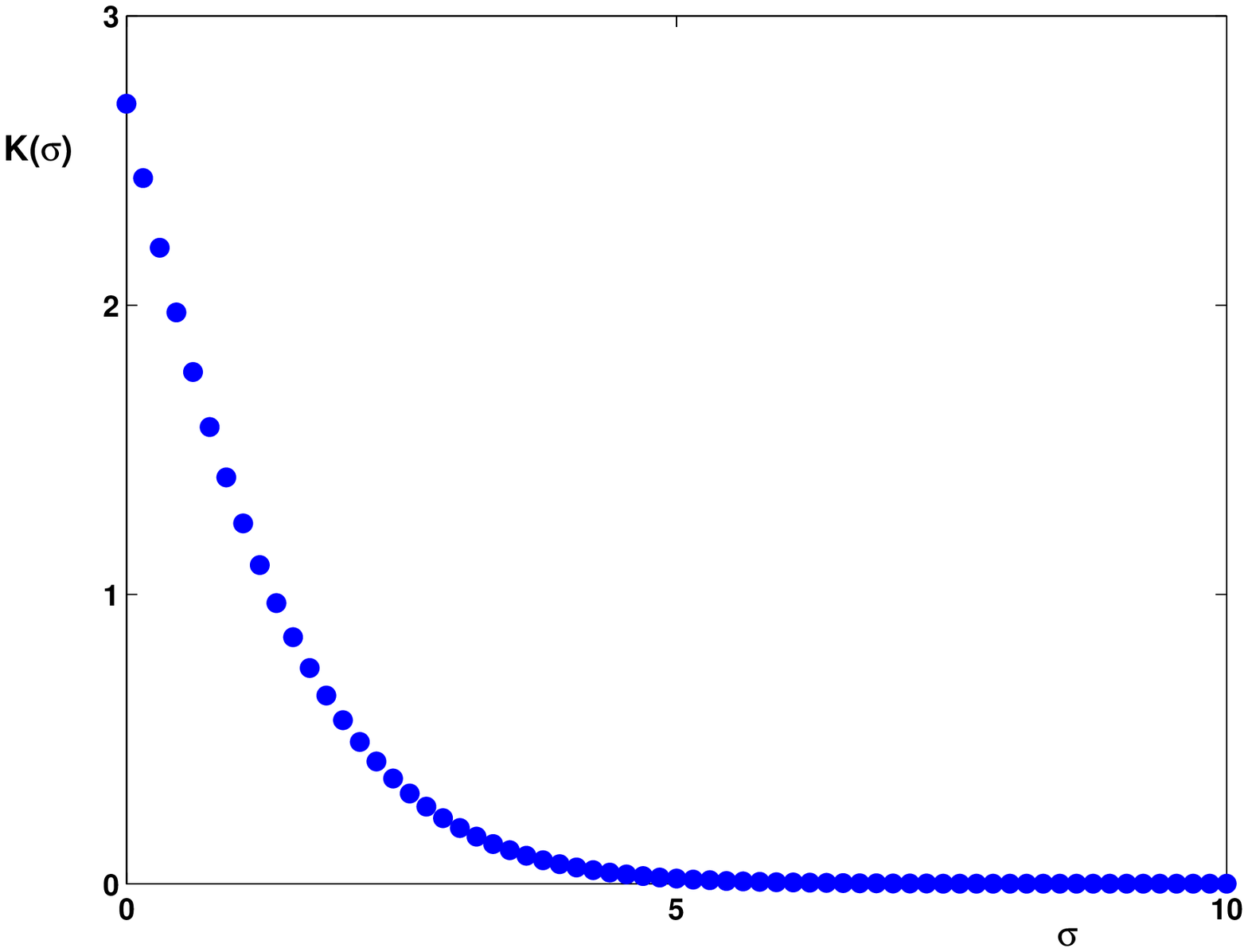} 
\includegraphics[height=5.5cm]{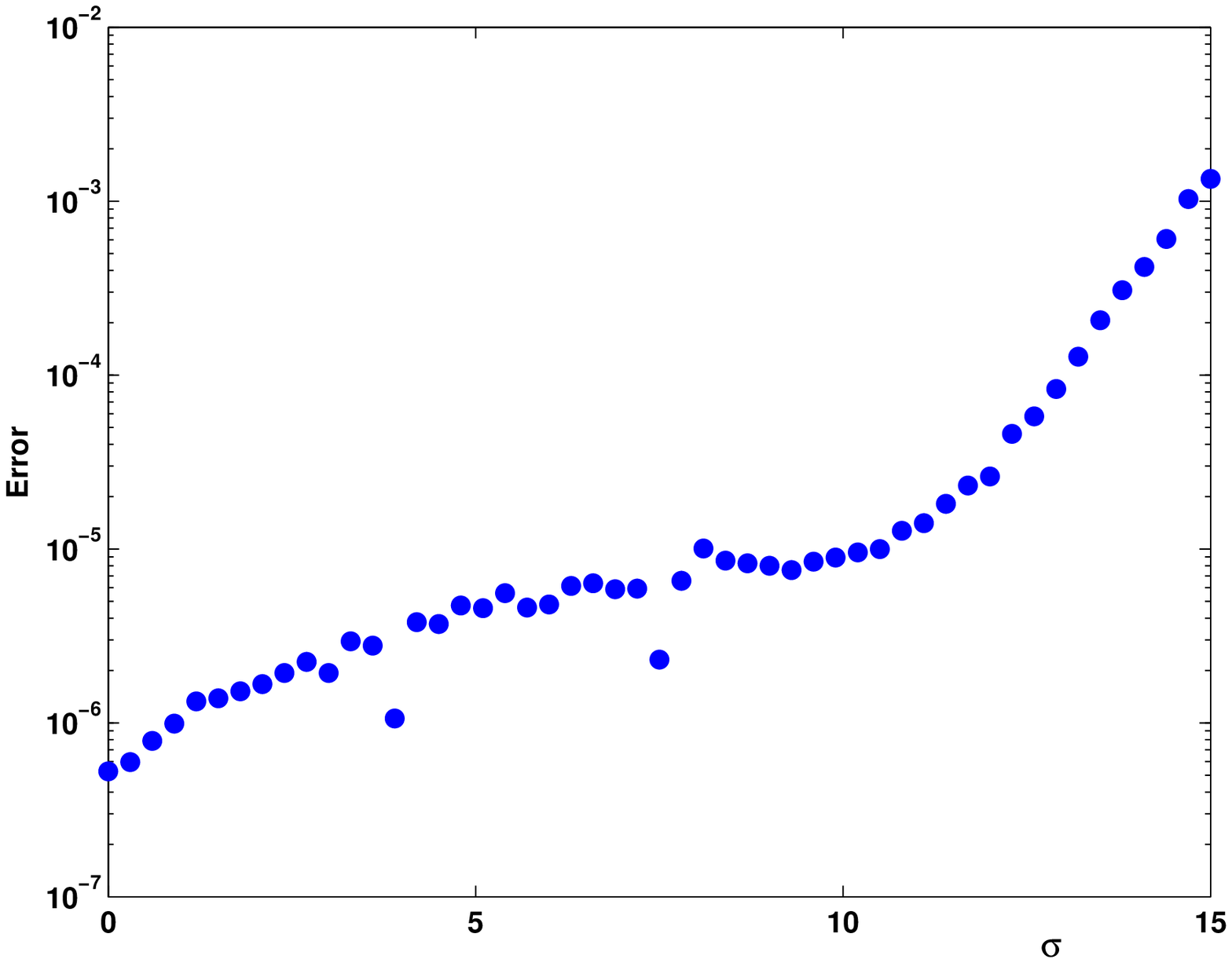}  
\end{center}
\caption{Left: The graph of the split function $K(\sigma)$ for $c_0 
= 0.00001$, $dt = 0.005$, and different values of $\sigma$ for kink 
solutions of the ODE (\ref{ODE}) with $F(\phi) = \sin(\phi)$. Right: 
relative error of numerical approximations of $K(\sigma)$ for $dt = 
0.005$ under variations of $c_0$.}
\end{figure}

Since the discrete sine--Gordon lattice has multiple equilibrium 
states at $\phi = \pi n$, $n \in \mathbb{Z}$, one can consider 
multiple monotonic kinks between several consequent equilibrium 
states \cite{CK00}. In particular, we shall consider the double 
kink, which represents a bound state between the heteroclinic orbit 
from $u_- = -\pi$ to $u_+ = \pi$ and the heteroclinic orbit from 
$u_+ = \pi$ to $\tilde{u}_+ = 3 \pi$. 

Double kinks were considered recently in a fourth-order differential 
model in \cite{CK00} and in the differential advance-delay equation 
in \cite{ACR03}. We note that the fourth-order model in \cite{CK00} 
is obtained in the continuous limit of the discrete Klein--Gordon 
lattice (\ref{Klein-Gordon}), such that the fourth-order derivative 
term is small compared to the second-order derivative term. However, 
numerical results of \cite{CK00} were derived in the case when both 
derivatives are comparable. Asymptotic and direct numerical results 
showed that the families of double kinks intersect the bifurcation 
point $(c,h) = (1,0)$, where our normal form 
(\ref{scalar-fourth-order}) is applicable. Thus, we give a rigorous 
explanation of existence of double kinks in the discrete sine-Gordon 
lattice from the theory of normal forms.

We modify the numerical algorithm by considering the initial-value
problem for the normalized fourth-order equation (\ref{ODE}) with
the initial values (\ref{initial-value}). Let $t_0$ be defined now
from the equation $\phi(t_0) = u_+$, and the split function
$K(\sigma)$ be defined by $K = \phi''(t_0)$. When $K = 0$, the
bounded solution $\phi_u(t)$ defines a double kink solution with the 
properties: $\lim_{t \to -\infty} \phi_u(t) = u_-$, $\lim_{t \to 
\infty} \phi_u(t) = \tilde{u}_+$ and $\phi_u(t - t_0) = 
-\phi_u(t_0-t) + 2 \pi$. The graph of $K$ versus $\sigma$ on Figure 
6 shows infinitely many zeros of $K(\sigma)$, which correspond to 
infinitely many double kink solutions (only four such solutions were 
reported in \cite{CK00}).

\begin{figure}[htbp]
\begin{center}
\includegraphics[height=8cm]{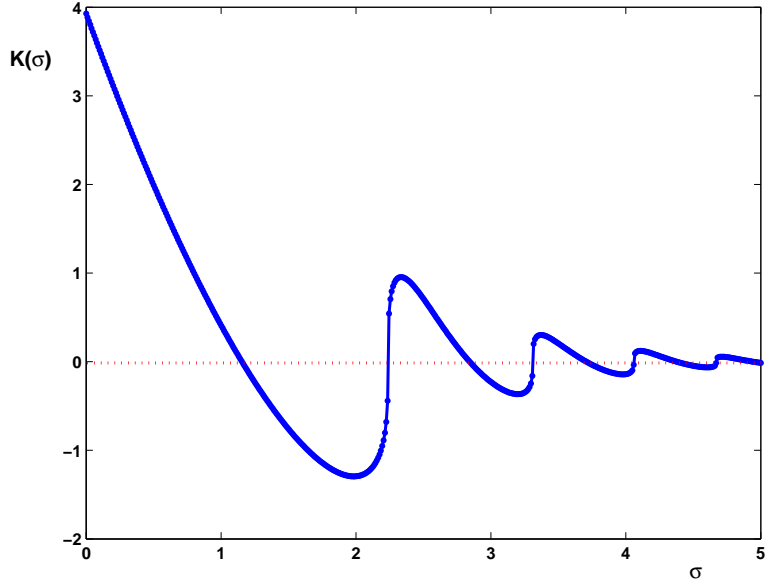}
\end{center}
\caption{The graph of the split function $K(\sigma)$ for double
kink solutions of the ODE (\ref{ODE}) with $F(\phi) =
\sin(\phi)$.}
\end{figure}       

In addition, we mention that the numerical method of computations of 
the split function can be tested for other heteroclinic orbits, such 
as triple monotonic kinks between $(-\pi,\pi)$, $(\pi,3\pi)$, 
$(3\pi,5\pi)$ and triple non-monotonic kinks between $(-\pi,\pi)$, 
$(\pi,-\pi)$, $(-\pi,\pi)$. It can also be used for construction of 
homoclinic orbits (see \cite{TP05}). The accuracy of computations of 
the split function is only limited to the accuracy of the ODE solver 
for the scalar fourth-order differential equation (\ref{ODE}).

\section{Persistence analysis of bounded solutions}

We consider the truncation of the system 
(\ref{center-manifold-Final}) in the form:
\begin{equation}
\dot{\mbox{\boldmath $\phi$}} = {\bf N}(\mbox{\boldmath $\phi$}),  
\label{reducedsyst}
\end{equation}%
where $\mbox{\boldmath $\phi$} = (\phi_1,\phi_2,\phi_3,\phi_4)^T$
and the vector field ${\bf N} : \mathbb{R}^4 \mapsto \mathbb{R}^4$ 
is given explicitly in the form
\begin{equation} 
{\bf N}(\mbox{\boldmath $\phi$}) = (\phi _{2},\phi _{3},\phi 
_{4},12\gamma \phi _{3}-12\tau F(\phi_1))^{T}.
\end{equation}
Let us use $t = \zeta_1$ in notations of this section. The truncated 
system (\ref{reducedsyst}) exhibits a number of useful properties:

\begin{itemize}
\item[P1] Since $F(\phi_1)$ is odd in $\phi_1 \in \mathbb{R}$, the vector 
field ${\bf N}(\mbox{\boldmath $\phi$})$ is odd in $\mbox{\boldmath 
$\phi$} \in \mathbb{R}^4$ 

\item[P2] The system is conservative, such that ${\rm div} \; {\bf N} = 
0$. 

\item[P3] The system is reversible under the 
reversibility symmetry ${\cal S}_0$ in 
(\ref{reversibility-reduced}). 

\item[P4] The set of fixed points of $\pm {\cal S}_0$ is two-dimensional.

\item[P5]  The system has two symmetric equilibrium points $\mathbf{u}_{\pm }=(u_{\pm 
},0,0,0)^{T}$. 

\item[P6] Since $F'(u_{\pm}) < 0$, the linearized operator $D {\bf N}({\bf u}_{\pm})$ 
has two simple eigenvalues on the imaginary axis and two simple real 
eigenvalues for any $\gamma \in \mathbb{R}$ and $\tau \neq 0$. 
\end{itemize}

We shall add a generic assumption on existence of a heteroclinic 
connection between periodic solutions near the non-zero equilibrium 
points. As it follows from numerical results in Section 4, we can 
not assume the existence of a true heteroclinic orbit between the 
two equilibrium points $\mathbf{u}_{\pm }$ since the true 
heteroclinic orbits are not supported by the nonlinear functions of 
the discrete $\phi^4$ and sine-Gordon equations. However, a 
numerical shooting method implies an existence of a one-parameter 
family of anti-reversible heteroclinic connections between periodic 
solutions. Indeed, an initial-value problem for an anti-reversible 
solutions $\mbox{\boldmath $\phi$}(t)$ of the truncated system 
(\ref{reducedsyst}) (such that ${\cal S}_0 \mbox{\boldmath 
$\phi$}(-t) = - \mbox{\boldmath $\phi$}(t)$) has two parameters, 
while there is only one constraint on the asymptotic behavior of 
$\mbox{\boldmath $\phi$}(t)$ as $t \to \infty$. An orthogonality 
constraint must be added to remove an eigendirection towards a 
one-dimensional unstable manifold associated to the real positive 
eigenvalue $\lambda_0$ in the linearization of ${\bf u}_+$. It is 
expected that a one-parameter family of solutions of the 
initial-value problem satisfies the given constraint. We will make a 
generic assumption on existence of the one-parameter family of 
anti-reversible solutions of the truncated system 
(\ref{reducedsyst}) and we intend to prove that such solutions 
persist for the full system (\ref{center-manifold-Final}).

\begin{Assumption}
\label{assumption-1} For a fixed set of parameter $(\gamma,\tau)$ 
and a fixed nonlinear function $F(\phi_1)$, there exists a smooth 
family of solutions $\mbox{\boldmath $\phi$}_{\alpha }(t)$, such 
that
\begin{eqnarray}
\dot{\mbox{\boldmath $\phi$}}_{\alpha } = {\bf N}(\mbox{\boldmath 
$\phi$}_{\alpha }), \label{phi_alpha} 
\end{eqnarray}
and
\begin{eqnarray}
{\cal S}_0 \mbox{\boldmath $\phi$}_{\alpha }(-t) = - \mbox{\boldmath 
$\phi$}_{\alpha }(t).
\end{eqnarray}
When $t\rightarrow \pm \infty$, the family $\mbox{\boldmath 
$\phi$}_{\alpha }(t)$ tends towards $\mathbf{u}_{\alpha }^{\pm 
}(t),$ where $\mathbf{u}_{\alpha }^{\pm }(t)$ are periodic solutions 
of the system (\ref{reducedsyst}) near the equilibrium points 
$\mathbf{u}_{\pm }$, such that
\begin{equation}
{\cal S}_{0}\mathbf{u}_{\alpha }^{-}(-t)=-\mathbf{u}_{\alpha 
}^{+}(t) \label{rev-sym-per}
\end{equation}
and
\begin{equation}
{\cal S}_{0}\mathbf{u}_{\alpha }^{-}(-t-\delta _{\alpha 
})=\mathbf{u}_{\alpha }^{-}(t-\delta _{\alpha }).  \label{rev-per}
\end{equation}%
We assume that $T_{\alpha }$ (period at infinity) and 
$\delta_{\alpha }$ (phase shift) are smooth functions of parameter 
$\alpha$. 
\end{Assumption}

Let us choose a particular heteroclinic connection $\mbox{\boldmath 
$\phi$}_{0}$. We define a linear operator $\mathcal{L}$ associated 
with the linearized system (\ref{reducedsyst}) near $\mbox{\boldmath 
$\phi$}_{0}(t)$:
\begin{equation}
\mathcal{L}(t) \mbox{\boldmath $\varphi$}(t) \doteq 
\dot{\mbox{\boldmath $\varphi$}}(t) - D {\bf N}(\mbox{\boldmath 
$\phi$}_{0}(t)) \mbox{\boldmath $\varphi$}(t). 
\label{linearized-system}
\end{equation}
If the periodic solutions ${\bf u}_0^{\pm}(t)$ have a sufficiently 
small radius, the spectrum of the linearized system 
(\ref{linearized-system}) on the periodic solutions ${\bf 
u}_0^{\pm}(t)$ is predicted from that at the equilibrium points 
${\bf u}^{\pm}$. Our aim is to study the invertibility of 
$\mathcal{L}$ in the space of continuous functions on $t \in 
\mathbb{R}$, which are asymptotically periodic as $t \to \pm 
\infty$. In what follows, we prove the following result, under an 
additional assumption that the spectrum of ${\cal L}$ is generic 
(free of bifurcation). The Fredholm Alternative Theorem would need 
to be used in the case when the assumption fails and the kernel of 
${\cal L}$ includes an anti-reversible decaying solution. 

\begin{Assumption}
\label{assumption-2} There is no anti-reversible solution of the 
homogeneous linearized system that decays to zero at both 
infinities, i.e. such that ${\cal L}(t) \mbox{\boldmath 
$\varphi$}(t) = {\bf 0}$, ${\cal S}_0 \mbox{\boldmath $\varphi$}(-t) 
= -\mbox{\boldmath $\varphi$}(t)$, and $\lim_{t \to \pm \infty} 
\mbox{\boldmath $\varphi$}(t) = {\bf 0}$. 
\end{Assumption}

\begin{Theorem} 
\label{theorem-persistence} Suppose Assumptions \ref{assumption-1} 
and \ref{assumption-2} hold. Let ${\bf F} \in 
C_{b}^{0}(\mathbb{R},\mathbb{R}^{4})$ be a function that is 
reversible with respect to ${\cal S}_0$ and tends towards a 
$T_{0}$-periodic function as $t \to \pm \infty$. Then, there exists 
a unique solution $\mbox{\boldmath $\varphi$} \in 
C_{b}^{0}(\mathbb{R},\mathbb{R}^{4})$ of the system
\begin{equation}
\label{proof-equation} \mathcal{L}(t) \mbox{\boldmath $\varphi$}(t) 
= {\bf F}(t),
\end{equation}
which is anti-reversible with respect to ${\cal S}_0$ and tends 
towards a $T_{0}$-periodic function as $t \to \pm \infty$. The 
linear map ${\bf F} \mapsto \mbox{\boldmath $\varphi$}$ is 
continuous in $C_{b}^{0}(\mathbb{R},\mathbb{R}^{4})$.
\end{Theorem}

\begin{Proof}
We shall first consider the fundamental matrix for the 
time-dependent linear system,
\begin{equation}
\label{proof-equation-1} \dot{\mbox{\boldmath $\varphi$}} = D {\bf 
N}(\mbox{\boldmath $\phi$}_{0}(t)) \mbox{\boldmath $\varphi$}.
\end{equation}
Due to properties of the truncated system (\ref{reducedsyst}), we 
have the following relationship,
\begin{equation}
{\cal S}_{0} D {\bf N}(\mbox{\boldmath $\phi$}_{0}(t)) 
\mbox{\boldmath $\varphi$} = - D {\bf N}(\mbox{\boldmath 
$\phi$}_{0}(-t)) {\cal S}_{0} \mbox{\boldmath $\varphi$}.
\end{equation}
By differentiating (\ref{phi_alpha}) with respect to $t$, we obtain 
the first eigenvector of ${\cal L}(t)$:
\begin{equation}
\mathcal{L}(t)\dot{\mbox{\boldmath $\phi$}}_{0}(t)=0.  
\label{phi_0dot}
\end{equation}%
By differentiation (\ref{phi_alpha}) with respect to $\alpha$, we 
can obtain another eigenvector of ${\cal L}(t)$. However, since the 
period depends on $\alpha$, the second eigenvector is non-periodic 
at infinity but contains a linearly growing term in $t$. So, let us 
define a new "time" $\tau$ as 
\begin{equation}
t=k_{\alpha } \tau, \qquad k_{\alpha }=T_{\alpha }/T_{0}, \qquad 
k_{0}=1.
\end{equation}%
Then $\widetilde{\mbox{\boldmath $\phi$}}_{\alpha }(\tau 
)=\mbox{\boldmath $\phi$}_{\alpha }(t)$ has the period $T_{0}$ as 
$\tau \rightarrow \pm \infty$ and satisfies
\begin{equation}
\frac{d}{d\tau }\widetilde{\mbox{\boldmath $\phi$}}_{\alpha } = 
k_{\alpha } {\bf N}(\widetilde{\mbox{\boldmath $\phi$}}_{\alpha }).
\end{equation}
Differentiating this identity with respect to $\alpha$ and denoting 
$\widehat{\mbox{\boldmath $\phi$}}_{0}=\partial _{\alpha 
}\widetilde{\mbox{\boldmath $\phi$}}_{\alpha }|_{\alpha =0}$ we 
obtain the equation,
\begin{equation}
\frac{d}{d\tau }\widehat{\mbox{\boldmath $\phi$}}_{0} = D {\bf 
N}(\mbox{\boldmath $\phi$}_{0})\widehat{ \mbox{\boldmath 
$\phi$}}_{0} + k_{1} {\bf N}(\mbox{\boldmath $\phi$}_{0}),
\end{equation}
where
\begin{equation}
k_{1}=\partial _{\alpha }k|_{\alpha =0}.
\end{equation}
This explicit computation shows that $\widehat{\mbox{\boldmath 
$\phi$}}_{0}(t)$ is the generalized eigenvector of the operator 
${\cal L}(t)$, 
\begin{equation}
\mathcal{L}(t)\widehat{\mbox{\boldmath $\phi$}}_{0}(t) = k_{1} 
\dot{\mbox{\boldmath $\phi$}}_{0}(t),  \label{Jordan block}
\end{equation}
which is bounded and $T_{0}$-periodic at infinity as the first 
eigenvector $\dot{\mbox{\boldmath $\phi$}}_{0}(t)$. Moreover, we 
note the symmetry properties:
\begin{eqnarray}
{\cal S}_{0} \dot{\mbox{\boldmath $\phi$}}_{0}(-t) &=& 
\dot{\mbox{\boldmath $\phi$}}_{0}(t), \\
{\cal S}_{0} \widehat{\mbox{\boldmath $\phi$}}_{0}(-t) &=& - 
\widehat{\mbox{\boldmath $\phi$}}_{0}(t).
\end{eqnarray}
The two eigenvectors $\dot{\mbox{\boldmath $\phi$}}_{0}(t)$ and 
$\widehat{\mbox{\boldmath $\phi$}}_{0}(t)$ correspond to the 
"neutral" directions associated with the families of periodic 
solutions $\mathbf{u}_{\alpha }^{\pm }(t)$ near the equilibria 
$\mathbf{u}_{\pm }$. Indeed, similarly we denote by 
$\mathbf{\dot{u}} _{0}^{\pm }$ and $\mathbf{\hat{u}}_{0}^{\pm }$ the 
$T_{0}-$ periodic solutions of the limiting system with periodic 
coefficients:
\begin{eqnarray*}
\frac{d}{dt}\mathbf{\dot{u}}_{0}^{\pm } &=& D {\bf N}(\mathbf{u}_{0}^{\pm })\mathbf{%
\dot{u}}_{0}^{\pm } \\
\frac{d}{dt}\mathbf{\hat{u}}_{0}^{\pm } &=& D {\bf N}(\mathbf{u}_{0}^{\pm })\mathbf{%
\hat{u}}_{0}^{\pm }+k_{1}\mathbf{\dot{u}}_{0}^{\pm }.
\end{eqnarray*}%
These two modes are the eigenvector and generalized eigenvector 
belonging to the double zero Floquet exponent associated with each 
periodic solution $\mathbf{u}_{0}^{\pm }(t)$. Since these periodic 
solutions are close enough to the equilibria $\mathbf{u}_{\pm },$ 
and thanks to perturbation theory, we know that there exist two real 
Floquet exponents $(\lambda _{0}^{\pm},-\lambda_0^{\pm})$ for each 
solution $\mathbf{u}_{0}^{\pm}$. Let us denote by $\mbox{\boldmath 
$\zeta$}_{\pm}(t)$ and $\mbox{\boldmath $\chi$}_{\pm}(t)$ the 
$T_{0}$-periodic eigenvectors (defined up to a factor) such that
\begin{eqnarray*}
\frac{d}{dt} \mbox{\boldmath $\zeta$}_{+} + \lambda _{0}^- 
\mbox{\boldmath $\zeta$}_{+} &=& D {\bf N}(\mathbf{u}_{0}^{\_}) \mbox{\boldmath $\zeta$}_{+} \\
\frac{d}{dt}\mbox{\boldmath $\zeta$}_{-}-\lambda _{0}^- 
\mbox{\boldmath $\zeta$}_{-} &=& D {\bf N}(\mathbf{u}_{0}^{\_}) \mbox{\boldmath $\zeta$}_{-} \\
\frac{d}{dt} \mbox{\boldmath $\chi$}_{+}+\lambda _{0}^+ 
\mbox{\boldmath $\chi$}_{+} &=& D {\bf 
N}(\mathbf{u}_{0}^{+}) \mbox{\boldmath $\chi$}_{+} \\
\frac{d}{dt} \mbox{\boldmath $\chi$}_{-} - \lambda _{0}^+ 
\mbox{\boldmath $\chi$}_{-} &=& D {\bf 
N}(\mathbf{u}_{0}^{+})\mbox{\boldmath $\chi$}_{-}.
\end{eqnarray*}%
We justify below that $\lambda _{0}^+ = \lambda_0^-$. Since the 
function $S_{0} \mbox{\boldmath $\zeta$}_{+}(-t)$ satisfies
\begin{equation}
\frac{d}{dt}S_{0} \mbox{\boldmath $\zeta$}_{+}(-t)- \lambda 
_{0}S_{0} \mbox{\boldmath $\zeta$}_{+}(-t) = D {\bf N}(S_{0} 
\mathbf{u}_{0}^{\_}(-t)) S_{0} \mbox{\boldmath $\zeta$}_{+}(-t)
\end{equation}%
and because of (\ref{rev-per}) we can choose $\mbox{\boldmath 
$\zeta$}_{-}$ such that
\begin{equation}
S_{0} \mbox{\boldmath $\zeta$}_{+}(-t-\delta _{0}) = \mbox{\boldmath 
$\zeta$}_{-}(t-\delta _{0}). \label{dzeta+ - dzeta-}
\end{equation}
Now, since we have (\ref{rev-sym-per}) and since $D {\bf 
N}(\mbox{\boldmath $\phi$}_0(t))$ is even in $\mbox{\boldmath 
$\phi$}_0(t)$, we deduce that
\begin{eqnarray*}
D {\bf N}(\mathbf{u}_{0}^{+}(t)) = D {\bf 
N}(-S_{0}(\mathbf{u}_{0}^{-}(-t)) = D {\bf 
N}(\mathbf{u}_{0}^{-}(t-2\delta _{0})),
\end{eqnarray*}
such that $\lambda _{0}^+ = \lambda_0^-$. Finally we can choose the 
vector functions $\mbox{\boldmath $\zeta$}_{\pm }$ and 
$\mbox{\boldmath $\chi$}_{\pm }$ such that
\begin{eqnarray}
\nonumber \mbox{\boldmath $\chi$}_{+}(t+\delta _{0}) &=& \mbox{\boldmath $\zeta$}_{+}(t-\delta _{0})  \\
\mbox{\boldmath $\chi$}_{-}(t+\delta _{0}) &=& \mbox{\boldmath 
$\zeta$}_{-}(t-\delta _{0})  \label{rel-xhi-dzeta} \\ \nonumber 
S_{0} \mbox{\boldmath $\chi$}_{+}(-t+\delta _{0}) &=& 
\mbox{\boldmath $\chi$} _{-}(t+\delta _{0}). 
\end{eqnarray}
Since $\mbox{\boldmath $\zeta$}_{\pm }(t)e^{\pm \lambda _{0}t}$ and 
$\mbox{\boldmath $\chi$}_{\pm }(t)e^{\pm \lambda _{0}t}$ are 
solutions of the limiting linear system, we introduce the unique 
solutions $\mbox{\boldmath $\phi$}_{1,2}(t)$ and $\mbox{\boldmath 
$\psi$}_{1,2}(t)$ of the system (\ref{proof-equation-1}) on $t \in  
\mathbb{R}$, such that 
\begin{eqnarray*}
&& \mbox{\boldmath $\phi$}_{1}(t) \sim 
\mbox{\boldmath $\zeta$}_+(t)e^{\lambda _{0}t} \quad \mbox{as} \quad  t\rightarrow -\infty \\
&& \mbox{\boldmath $\phi$}_{2}(t) \sim 
\mbox{\boldmath $\zeta$}_{-}(t)e^{-\lambda _{0}t}  \quad \mbox{as} \quad  t\rightarrow -\infty  \\
&& \mbox{\boldmath $\psi$}_{2}(t) \sim \mbox{\boldmath $\chi$}_{+}(t)e^{\lambda _{0}t}  
\quad \mbox{as} \quad  t\rightarrow \infty  \\
&& \mbox{\boldmath $\psi$}_{1}(t) \sim \mbox{\boldmath 
$\chi$}_{-}(t)e^{-\lambda _{0}t}  \quad \mbox{as} \quad  
t\rightarrow \infty.
\end{eqnarray*}%
Now we have%
\begin{eqnarray*}
\frac{d}{dt}S_{0} \mbox{\boldmath $\phi$}_{1}(-t) = D {\bf N}(S_{0} 
\mbox{\boldmath $\phi$}_{0}(-t))S_{0} \mbox{\boldmath 
$\phi$}_{1}(-t) = D {\bf N}(\mbox{\boldmath 
$\phi$}_{0}(t))S_{0}\mbox{\boldmath $\phi$}_{1}(-t)
\end{eqnarray*}%
and%
\begin{eqnarray*}
S_{0} \mbox{\boldmath $\phi$}_{1}(-t) \sim S_{0}\mbox{\boldmath 
$\zeta$}_{+}(-t) e^{-\lambda_{0}t}   \quad \mbox{as} \quad  
t\rightarrow \infty,
\end{eqnarray*}
where
\begin{eqnarray*}
S_{0} \mbox{\boldmath $\zeta$}_{+}(-t)e^{-\lambda_{0}t} = 
\mbox{\boldmath $\zeta$}_{-}(t-2\delta _{0})e^{-\lambda _{0}t} = 
\mbox{\boldmath $\chi$}_{-}(t)e^{-\lambda _{0}t}.
\end{eqnarray*}%
Hence%
\begin{equation}
S_{0} \mbox{\boldmath $\phi$}_{1}(-t)=\mbox{\boldmath 
$\psi$}_{1}(t),  \label{Sphi_1}
\end{equation}%
and, in the same way%
\begin{equation}
S_{0}\mbox{\boldmath $\phi$}_{2}(-t)=\mbox{\boldmath $\psi$}_{2}(t).  
\label{Sphi_2}
\end{equation}
The two pairs of linearly independent solutions $( \mbox{\boldmath 
$\phi$}_{1},\mbox{\boldmath $\phi$}_{2})$ and $(\mbox{\boldmath 
$\psi$}_{1},\mbox{\boldmath $\psi$}_{2})$ cannot contain components 
in $\dot{\mbox{\boldmath $\phi$}}_{0}$ and $\widehat{\mbox{\boldmath 
$\phi$}}_{0}$ thanks to the behavior at infinity, hence there exist 
constants $(a,b,c,d)$ such that 
\begin{eqnarray}
\nonumber \mbox{\boldmath $\psi$}_{1} &=& a \mbox{\boldmath $\phi$}_{1} + 
b \mbox{\boldmath $\phi$}_{2} \\
\label{rel-psi-phi-0} \mbox{\boldmath $\psi$}_{2} &=& c 
\mbox{\boldmath $\phi$}_{1} + d \mbox{\boldmath $\phi$}_{2}.
\end{eqnarray}%
Iterating once the relations (\ref{Sphi_1}) and (\ref{Sphi_2}), we 
find that 
$$
\left( \begin{array}{cc} a & b \\ c & d
\end{array} \right) ^{2} = \left( \begin{array}{cc}
a^2 + bc & b(a+d) \\ c(a+d) & d^2 + bc \end{array} \right) = \left( 
\begin{array}{cc} 1 & 0 \\ 0 & 1 \end{array} \right).
$$
If either $b$ or $c$ is non-zero, then $d = -a$ and $a^2 + bc = 1$. 
If $b = c = 0$, then $a^2 = d^2 = 1$. The choice $a = - d = \pm 1$ 
is included in the previous case $d = -a$. The choice $a = d = \pm 
1$ is impossible since it would lead to more than two independent 
vectors in the two-dimensional set of fixed points of the operators 
$\pm {\cal S}_0$ at $t = 0$. Thus, we parameterize: $a = -d = 
\alpha$, $b = \beta$, and $c = \gamma$, where $\alpha^2 + \beta 
\gamma = 1$ and rewrite the relationship (\ref{rel-psi-phi-0}) in 
the form:
\begin{eqnarray}
\nonumber \mbox{\boldmath $\psi$}_{1} &=& \alpha \mbox{\boldmath $\phi$}_{1} 
+ \beta \mbox{\boldmath $\phi$}_{2} \\
\label{rel-psi-phi} \mbox{\boldmath $\psi$}_{2} &=& \gamma 
\mbox{\boldmath $\phi$}_{1} - \alpha \mbox{\boldmath $\phi$}_{2}
\end{eqnarray}%
with the inverse relationship:
\begin{eqnarray}
\nonumber \mbox{\boldmath $\phi$}_{1} &=&\alpha \mbox{\boldmath 
$\psi$}_{1} + \beta \mbox{\boldmath $\psi$}_{2}  \\
\label{rel-phi-psi} \mbox{\boldmath $\phi$}_{2} &=&\gamma 
\mbox{\boldmath $\psi$}_{1}-\alpha \mbox{\boldmath $\psi$}_{2}.  
\end{eqnarray}
We now wish to solve the linear equation (\ref{proof-equation}), 
where the right-hand-side vector $\mathbf{F}(t)$ is asymptotically 
$T_0$-periodic and reversible with respect to ${\cal S}_0$, i.e. 
such that
\begin{equation}
S_{0}\mathbf{F}(-t)=\mathbf{F}(t).  \label{revF}
\end{equation}%
By construction, we have%
\begin{eqnarray*}
\mathbf{F} &=&f_{0} \dot{\mbox{\boldmath $\phi$}}_{0} + f_{1} 
\widehat{\mbox{\boldmath $\phi$}}_{0} + g_{1} \mbox{\boldmath 
$\phi$}_1 + g_{2} \mbox{\boldmath $\phi$}_{2} \\
\mbox{\boldmath $\varphi$} &=&\alpha _{0}\dot{\mbox{\boldmath 
$\phi$}}_{0}+ \alpha_{1} \widehat{\mbox{\boldmath $\phi$}}_{0}+\beta 
_{1} \mbox{\boldmath $\phi$}_{1}+\beta_{2} \mbox{\boldmath 
$\phi$}_{2},
\end{eqnarray*}%
where the functions $g_{1}(t)e^{-\lambda _{0}|t|}$ and 
$g_{2}(t)e^{\lambda _{0}|t|}$ are bounded as $t\rightarrow -\infty$ 
and asymptotically $T_{0}$-periodic. Let $\mathcal{M}(t) = 
[\dot{\mbox{\boldmath $\phi$}}_{0},\hat{\mbox{\boldmath 
$\phi$}}_{0},\mbox{\boldmath $\phi$}_{1},\mbox{\boldmath 
$\phi$}_{2}]$, such that the coordinates of $\mathbf{F}$ and 
$\mbox{\boldmath $\varphi$}$ in the new basis are
\begin{eqnarray}
\widetilde{\mathbf{F}}(t) &=&\mathcal{M}(t)^{-1}\mathbf{F}%
(t)=(f_{0},f_{1},g_{1},g_{2})^T  \label{compF,Phi} \\
\widetilde{\mbox{\boldmath $\varphi$}}(t) &=& \mathcal{M}(t)^{-1} 
\mbox{\boldmath $\varphi$}(t)=(\alpha _{0},\alpha _{1},\beta 
_{1},\beta _{2})^T.
\end{eqnarray}
The system (\ref{proof-equation}) reduces to the simple system
\begin{eqnarray*}
\frac{d\alpha _{0}}{dt} &=&-k_{1}\alpha _{1}+f_{0} \\
\frac{d\alpha _{1}}{dt} &=&f_{1} \\
\frac{d\beta _{1}}{dt} &=&g_{1} \\
\frac{d\beta _{2}}{dt} &=&g_{2}.
\end{eqnarray*}%
Since the system (\ref{reducedsyst}) is conservative (that is ${\rm 
tr}\left( D {\bf N}(\mbox{\boldmath $\phi$}_{0})\right)=0$), the 
Liouville's theorem implies that 
\begin{equation}
\det \mathcal{M}(t)=\det \mathcal{M}(0)\neq 0.
\end{equation}
Solving the uncoupled system for $\beta_{1,2}(t)$, we find that 
\begin{eqnarray*}
\beta _{2}(t) &=&\int_{-\infty }^{t}g_{2}(s)ds, \\
\beta _{1}(t) &=&\int_{0}^{t}g_{1}(s)ds+\beta _{10},
\end{eqnarray*}
such that the functions $\beta _{1}(t)e^{-\lambda _{0}|t|}$ and 
$\beta _{2}(t)e^{\lambda _{0}|t|}$ are bounded as $t\rightarrow 
-\infty$ and asymptotically $T_{0}$-periodic for any constant 
$\beta_{10}$. The same property holds for $\mbox{\boldmath 
$\varphi$}(t)$ as $t\rightarrow -\infty$. Now, we consider the 
opposite limit $t \to +\infty$. Thanks to the relations 
(\ref{rel-psi-phi}), we notice that the functions $(\alpha 
g_{1}+\gamma g_{2})(t)e^{-\lambda _{0}t}$ and $(\beta g_{1}-\alpha 
g_{2})(t)e^{\lambda _{0}t}$ are bounded as $t\rightarrow \infty$. 
Moreover they are asymptotically $T_{0}$-periodic if $\mathbf{F}(t)$ 
has this property. Now we can write the function $\beta _{1} 
\mbox{\boldmath $\phi$}_{1} + \beta _{2} \mbox{\boldmath 
$\phi$}_{2}$ as
\begin{eqnarray*}
\beta _{1} \mbox{\boldmath $\phi$}_{1} + \beta _{2} \mbox{\boldmath 
$\phi$}_{2} &=& \mbox{\boldmath $\psi$}_{1}(t)\left( 
\int_{0}^{t}(\alpha g_{1}(s) + \gamma g_{2}(s))ds + 
\int_{-\infty}^{0} \gamma g_{2}(s)ds+\alpha \beta_{10}\right) + \\
&&+ \mbox{\boldmath $\psi$}_{2}(t)\left( \int_{0}^{t}(\beta 
g_{1}(s)-\alpha g_{2}(s))ds+\int_{-\infty }^{0}-\alpha 
g_{2}(s)ds+\beta \beta _{10}\right) 
\end{eqnarray*}%
and it is easy to check the boundedness of the function $\beta_{1} 
\mbox{\boldmath $\phi$}_{1} + \beta _{2} \mbox{\boldmath 
$\phi$}_{2}$ as $t \to \infty$, provided that
\begin{equation}
\label{constraint-hold} \int_{0}^{\infty }(\beta g_{1}(s)-\alpha 
g_{2}(s))ds-\int_{-\infty }^{0}\alpha g_{2}(s)ds+\beta \beta 
_{10}=0.
\end{equation}
If $\beta \neq 0$, this equation determines the value of 
$\beta_{10}$. If $\beta = 0$ and $\alpha = +1$, the value of 
$\beta_{10}$ is not determined but equation (\ref{constraint-hold}) 
is satisfied identically since $g_2(t)$ is odd, thanks to the 
identity (\ref{symg_j}) below. The case $\beta = 0$ and $\alpha = 
-1$ is excluded by Assumption \ref{assumption-2}. Thus, the function 
$\beta _{1} \mbox{\boldmath $\phi$}_{1} + \beta_{2} \mbox{\boldmath 
$\phi$}_{2}$ is bounded as $t \to +\infty$ and asymptotically 
$T_{0}$-periodic. It remains to show that the reversibility 
constraint (\ref{revF}) implies the anti-reversibility constraint 
for the solution $\mbox{\boldmath $\varphi$}(t)$:
\begin{equation}
\mathcal{S}_{0} \mbox{\boldmath $\varphi$}(-t) = - \mbox{\boldmath 
$\varphi$}(t).
\end{equation}
By construction, the functions $f_{0}(t)$ and $f_{1}(t)$ are even 
and odd respectively. Since they are defined by a multiplication of 
an exponentially growing function with an exponentially decaying one 
of the same rate and with periodic factors, these functions are 
bounded and asymptotically $T_{0}$-periodic. Solving the uncoupled 
system for $\alpha_0(t)$ and $\alpha_1(t)$, we have 
\begin{eqnarray*}
\alpha _{1}(t) &=&\alpha _{1}(0)+\int_{0}^{t}f_{1}(s)ds, \\
\alpha _{0}(t) &=&\int_{0}^{t} \left[ f_{0}(s)-k_{1}\alpha _{1}(s) 
\right] ds.
\end{eqnarray*}
Since $f_{0}(t)$ is even and $f_{1}(t)$ is odd, it is clear that 
$\alpha _{1}(t)$ is even and $\alpha_{0}(t)$ is odd. Moreover, since 
$f_{1}(t)$ is odd, it has asymptotically a zero average at infinity. 
Therefore, the function $\alpha _{1}(t)$ is bounded and 
asymptotically $T_{0}$-periodic at infinity. Now we can find a 
unique $\alpha _{1}(0)$ such that the even function 
$f_{0}(s)-k_{1}\alpha _{1}(s)$ which is asymptotically 
$T_{0}$-periodic at infinity, has in addition a zero average at 
infinity. Then this gives a bounded function $\alpha _{0}(t)$ which 
is asymptotically $T_{0}$-periodic at infinity. Now for the two last 
components of $\mathbf{F}(t)$ we have
\begin{equation}
g_{1}(-t) \mbox{\boldmath $\psi$}_{1}(t)+g_{2}(-t)  \mbox{\boldmath 
$\psi$}_{2}(t) - g_{1}(t) \mbox{\boldmath $\phi$}_{1}(t) - g_{2}(t) 
\mbox{\boldmath $\phi$}_{2}(t)=0,
\end{equation}
which corresponds to the identities
\begin{eqnarray}
g_{1}(t)-\alpha g_{1}(-t)-\gamma g_{2}(-t) &=&0, \nonumber  \\
g_{2}(t)-\beta g_{1}(-t)+\alpha g_{2}(-t) &=&0. \label{symg_j}
\end{eqnarray}
We need to show that this implies%
\begin{eqnarray*}
\beta _{1}(t)+\alpha \beta _{1}(-t)+\gamma \beta _{2}(-t) &=&0, \\
\beta _{2}(t)+\beta \beta _{1}(-t)-\alpha \beta _{2}(-t) &=&0.
\end{eqnarray*}%
The identities (\ref{symg_j}) lead to%
\begin{equation}
\beta g_{1}(s)-\alpha g_{2}(s)=g_{2}(-s)
\end{equation}%
hence the constant $\beta _{10}$ takes the form%
\begin{equation}
\label{help-eq-111} \beta \beta _{10}=(\alpha -1)\int_{-\infty 
}^{0}g_{2}(s)ds.
\end{equation}
We observe that the identity (\ref{help-eq-111}) is satisfied when 
$\beta = 0$ and $\alpha = 1$. Now we check that
\begin{eqnarray*}
&& \beta _{1}(t)+\alpha \beta _{1}(-t)+\gamma \beta _{2}(-t)=\\
&=&\int_{0}^{t}\{g_{1}(s)-\alpha g_{1}(-s)\}ds+\gamma 
\int_{t}^{\infty
}g_{2}(-s)ds+(1+\alpha )\beta _{10} \\
&=&\gamma \int_{0}^{\infty }g_{2}(-s)ds+(1+\alpha )\beta _{10} 
=(1+\alpha )\beta _{10}+\gamma \int_{-\infty }^{0}g_{2}(s)ds.
\end{eqnarray*} 
If $\beta = 0$ and $\alpha = 1$, the expression above is zero under 
a special choice of $\beta_{10}$. If $\beta \neq 0$ and $\gamma = 
(1-\alpha^2)/\beta$, it is zero because
\begin{eqnarray*}
(1+\alpha )\beta _{10}+\gamma \int_{-\infty }^{0}g_{2}(s)ds = 
\frac{1+\alpha }{\beta }\left( \beta \beta _{10}+(1-\alpha 
)\int_{-\infty }^{0}g_{2}(s)ds\right) =0.
\end{eqnarray*}%
Similarly,
\begin{eqnarray*}
&& \beta _{2}(t)+\beta \beta _{1}(-t)-\alpha \beta _{2}(-t)= \\
&=&\int_{-\infty }^{t}g_{2}(s)ds+\beta \beta _{10}-\alpha 
\int_{-\infty}^{0} g_{2}(s) ds + \int_{0}^{t} 
\left[ \alpha g_{2}(-s)-\beta g_{1}(-s) \right] ds \\
&=&\beta \beta _{10}+(1-\alpha )\int_{-\infty }^{0}g_{2}(s)ds=0.
\end{eqnarray*}%
Hence the statement is proved in all cases except the case $\beta = 
0$ and $\alpha = -1$, which is excluded by Assumption 
\ref{assumption-2}. 
\end{Proof}

Using analysis of the truncated vector system (\ref{reducedsyst}), 
we can now prove persistence of bounded solutions of the full 
system. Since there exists a smooth mapping $\mbox{\boldmath 
$\psi$}(\mbox{\boldmath $\phi$},\varepsilon)$ in the hyperbolic 
subspace (see Theorem \ref{theorem-center-manifold}), we rewrite the 
system (\ref{center-manifold-Final}) in the form,
\begin{equation}
\dot{\mbox{\boldmath $\phi$}} = {\bf N}(\mbox{\boldmath $\phi$}) + 
\sqrt{\varepsilon } {\bf R}(\mbox{\boldmath $\phi$},\mbox{\boldmath 
$\psi$}(\mbox{\boldmath $\phi$},\varepsilon)),  \label{perturbed 
syst}
\end{equation}
where $\mbox{\boldmath $\psi$}(\mbox{\boldmath $\phi$},\varepsilon ) 
= {\rm O}(\sqrt{\varepsilon })$. The perturbation term is non-local, 
which forces us to use a non-geometric proof for \ showing the 
persistence of heteroclinic connections to periodic solutions near 
equilibria $\mathbf{u}_{\pm }$. Indeed, let us decompose 
$\mbox{\boldmath $\phi$}$ as
\begin{equation}
\mbox{\boldmath $\phi$} = \mbox{\boldmath $\phi$}_{0} + 
\mbox{\boldmath $\varphi$}
\end{equation}
and define
\begin{equation}
{\bf R}_{1}(\mbox{\boldmath $\phi$}_{0},\mbox{\boldmath $\varphi$}) 
= {\bf N}(\mbox{\boldmath $\phi$}_{0}+ \mbox{\boldmath $\varphi$}) - 
{\bf N}(\mbox{\boldmath $\phi$}_{0}) - D {\bf N}(\mbox{\boldmath 
$\phi$}_{0}) \mbox{\boldmath $\varphi$} = {\rm O}(||\mbox{\boldmath 
$\varphi$} ||^{2}),
\end{equation}
then $\mbox{\boldmath $\varphi$}$ satisfies
\begin{equation}
\mathcal{L}(t) \mbox{\boldmath $\varphi$} = {\bf 
R}_{1}(\mbox{\boldmath $\phi$}_{0},\mbox{\boldmath $\varphi$}) + 
\sqrt{\varepsilon } {\bf R}(\mbox{\boldmath 
$\phi$}_{0}+\mbox{\boldmath $\varphi$},\mbox{\boldmath 
$\psi$}(\mbox{\boldmath $\phi$}_{0}+\mbox{\boldmath 
$\varphi$},\varepsilon )).  \label{equperturbphi}
\end{equation}%
Using Theorem \ref{theorem-persistence} and the smooth dependence of 
${\bf R}$ on $(\mbox{\boldmath $\phi$},\mbox{\boldmath $\psi$})$, we 
can then use the implicit function theorem near $\mbox{\boldmath 
$\varphi$} = {\bf 0}$ and find a unique solution $\mbox{\boldmath 
$\varphi$} \in C_{b}^{0}(\mathbb{R},\mathbb{R}^{4})$ for small 
enough $\varepsilon$. Moreover, this solution $\mbox{\boldmath 
$\varphi$}(t)$ is a $T_{0}$-periodic function at infinity, which 
magnitude is bounded by a constant of the order of ${\rm 
O}(\sqrt{\varepsilon})$. Using these observations, we assert the 
following theorem.

\begin{Theorem}
\label{theorem-bounded-solutions} Let Assumptions \ref{assumption-1} 
and \ref{assumption-2} be satisfied and $\mbox{\boldmath 
$\phi$}_{0}(t)$ be a solution of the truncated system 
(\ref{reducedsyst}), which is heteroclinic to small $T_{0}$-periodic 
solutions ${\bf u}_{\pm}(t)$. Then, for small enough $\varepsilon$, 
there exists a unique solution $\mbox{\boldmath $\phi$}(t)$ of the 
perturbed system (\ref{perturbed syst}), which is 
$\sqrt{\varepsilon}$-close to $\mbox{\boldmath $\phi$}_{0}(t)$ and 
is heteroclinic to small $T_{0}$-periodic solutions of the system 
(\ref{perturbed syst}).
\end{Theorem}

We note that uniqueness is a remarkable property of the perturbed 
system (\ref{perturbed syst}). It is supported here due to the 
choice that the function $\mbox{\boldmath $\phi$}(t)$ has the same 
asymptotic period at infinity as the function $\mbox{\boldmath 
$\phi$}_{0}(t)$.

\section{Conclusion}

We have derived a scalar normal form equation for bifurcations of 
heteroclinic orbits in the discrete Klein-Gordon equation. Existence 
of the center manifold and persistence of bounded solutions of the 
normal form equation are proved with rigorous analysis. Bounded 
solutions may include heteroclinic orbits between periodic 
perturbations at the equilibrium states and true heteroclinic orbits 
between equilibrium states. Our numerical results indicate that no 
true heteroclinic orbits between equilibrium states exist for two 
important models, the discrete $\phi^4$ and sine--Gordon equations. 
Double and multiple kinks (between several equilibrium states) may 
exist in the case of the discrete sine-Gordon equation. 

We have studied persistence of heteroclinic anti-reversible 
connections between periodic solutions near the equilibrium states, 
but we have not addressed persistence of true heteroclinic orbits 
between the equilibrium states. It is intuitively clear that 
non-existence of a true heteroclinic orbit in the truncated normal 
form must imply the non-existence of such an orbit in the 
untruncated system at least for sufficient small $\varepsilon$. 
Conversely, existence of a true heteroclinic orbit in the truncated 
system should imply a continuation of the family for the full system 
on the two-parameter plane $(c,h)$. Rigorous analysis of such 
continuations is based on the implicit function argument for a split 
function, which becomes technically complicated since our full 
system involves the advance-delay operators, where the resolvent 
estimates are not simple. We pose this problem as an open question 
for further studies.

\appendix
\section{Truncation of the vector normal form from \cite{I95}}

The eigenvectors for the quadruple zero eigenvalues satisfy the
reversibility symmetry relations
(\ref{reversibility-eigenvectors}). The general reversible normal
form for this case was derived in \cite{I95}. We will show here
that this normal form can be reduced to the fourth-order equation
(\ref{scalar-fourth-order}) under an appropriate scaling. The
general reversible normal form is explicitly written in the vector
form (see \cite{I95} for details):
\begin{eqnarray*}
\xi_1' & = & \xi_2 + u_3 P_1(u_1,u_2,u_4), \\
\xi_2' & = & \xi_3 + q_3 P_1(u_1,u_2,u_4) + u_1 P_2(u_1,u_2,u_4) +
u_2 P_3(u_2,u_4), \\
\xi_3' & = & \xi_4 + r_3 P_1(u_1,u_2,u_4) + q_1 P_2(u_1,u_2,u_4) +
q_2 P_3(u_2,u_4) + u_3 P_4(u_1,u_2,u_4), \\
\xi_4' & = & s_3 P_1(u_1,u_2,u_4) + r_1 P_2(u_1,u_2,u_4) + r_2
P_3(u_2,u_4) + q_3 P_4(u_1,u_2,u_4) + P_5(u_1,u_2,u_4),
\end{eqnarray*}
where $P_{1,2,3,4,5}$ are polynomials in variables:
\begin{eqnarray*}
u_1 & = & \xi_1, \quad u_2 = \xi_2^2 - 2 \xi_1 \xi_3, \quad
u_3 = \xi_2^3 - 3 \xi_1 \xi_2 \xi_3 + 3 \xi_1^2 \xi_4, \\
u_4 & = & 3 \xi_2^2 \xi_3^2 - 6 \xi_2^3 \xi_4 - 8 \xi_1 \xi_3^3 +
18 \xi_1 \xi_2 \xi_3 \xi_4 - 9 \xi_1^2 \xi_4^2,\\
q_1 & = & \xi_2, \quad q_2 = -3\xi_1 \xi_4 + \xi_2 \xi_3,  \quad
q_3 = 3 \xi_1 \xi_2 \xi_4 - 4 \xi_1 \xi_3^2 + \xi_2^2 \xi_3, \\
r_1 & = & \xi_3, \quad r_2 = -3 \xi_2 \xi_4 + 2 \xi_3^2, \quad r_3
= - 3 \xi_1 \xi_3 \xi_4 + 3 \xi_2^2 \xi_4 - \xi_2 \xi_3^2,
\end{eqnarray*}
and $s_3 = 3 \xi_2 \xi_3 \xi_4 - \frac{4}{3} \xi_3^3 - 3 \xi_1
\xi_4^2$. If the zero equilibrium point persists in the vector 
normal form, then $P_5(0,0,0) = 0$. The linear part of the vector 
normal form corresponds to the truncation of $P_{1,2,3,4,5}$ as 
follows:
$$
P_1 = 0, \quad P_2 = \mu_1, \quad P_3 = P_4 = 0, \quad P_5 = \mu_2
u_1,
$$
where $\mu_1$ and $\mu_2$ are two parameters of the linear system.
The two parameters $(\mu_1,\mu_2)$ must recover the Taylor series
expansion (\ref{quadric-curve}) of the full dispersion relation
(\ref{dispersion-relation}) in scaled variables. Since the
dispersion relation from the linear part of the vector normal form
is
$$
\Lambda^4 - 3 \mu_1 \Lambda^2 + \mu_1^2 - \mu_2 = 0,
$$
the comparison with the Taylor series expansion
(\ref{quadric-curve}) shows that
\begin{equation}
\mu_1 = 4 \varepsilon \gamma + {\rm O}(\varepsilon^2), \qquad
\mu_2 = \varepsilon^2 (16 \gamma^2 - 12 \tau) + {\rm
O}(\varepsilon^3).
\end{equation}
The nonlinear part in the dynamical system
(\ref{dynamical-system-bifurcation}) is scaled as $\varepsilon^2$
near the bifurcation point (\ref{perturbed-point-bifurcation}), such 
that  
$$
P_1 = \varepsilon^2 \tilde{P}_1, \quad P_2 = \mu_1 + \varepsilon^2 
\tilde{P}_2, \quad P_{3,4} = \varepsilon^2 \tilde{P}_{3,4}, \quad 
P_4 = \varepsilon^2 \tilde{P}_4, \quad P_5 = \mu_2 u_1 + 
\varepsilon^2 \tilde{P}_5,
$$
where tilded functions $\tilde{P}_j$ depend on the same variables as 
the original functions $P_j$. Let $\zeta_1 = \sqrt{\varepsilon} 
\zeta$ and use the scaling transformation:
$$
\xi_1 = \phi_1(\zeta_1), \quad \xi_2 = \sqrt{\varepsilon} 
\phi_2(\zeta_1), \quad \xi_3 = \varepsilon \phi_3(\zeta_1), \quad 
\xi_4 = \sqrt{\varepsilon^3} \phi_4(\zeta_1).
$$
The scaled version of the vector normal form is written explicitly 
as follows:
\begin{eqnarray*}
\phi_1' & = & \phi_2 + \varepsilon^3 \tilde{u}_3 \tilde{P}_1, \\
\phi_2' & = & \phi_3 + 4 \gamma \phi_1 + \varepsilon \phi_1
\tilde{P}_2 + \varepsilon^2 \tilde{u}_2 \tilde{P}_3 + \varepsilon^3 \tilde{q}_3 \tilde{P}_1, \\
\phi_3' & = & \phi_4 + 4 \gamma \phi_2 + \varepsilon \phi_2 
\tilde{P}_2 + \varepsilon^2 ( \tilde{q}_2 \tilde{P}_3 + \tilde{u}_3 
\tilde{P}_4) + \varepsilon^3 \tilde{r}_3 \tilde{P}_1, \\
\phi_4' & = & 4 \gamma \phi_3 + (16 \gamma^2 - 12 \tau) \phi_1 + 
\tilde{P}_5 + \varepsilon \phi_3 \tilde{P}_2 + \varepsilon^2 ( 
\tilde{r}_2 \tilde{P}_3 + \tilde{q}_3 \tilde{P}_4) + \varepsilon^3 
\tilde{s}_3 \tilde{P}_1,
\end{eqnarray*}
where tilded variables stand for the scaled versions of the original 
variables. Let $G(\phi_1) = \tilde{P}_5(\phi_1,0,0)$. The truncated 
version of the vector normal form (with $\varepsilon \equiv 0$) is 
equivalent to the scalar fourth-order equation:
\begin{equation}
\label{scalar-fourth-order-equation} \phi^{({\rm iv})} - 12 \gamma
\phi'' + 12 \tau \phi = G(\phi),
\end{equation}
where $\phi \equiv \phi_1$. The comparison with the scalar 
fourth-order equation (\ref{scalar-fourth-order}) shows that 
$G(\phi) = -12 \tau \left( F(\phi) - \phi \right)$, where $F(\phi)$ 
is the nonlinearity of the Klein--Gordon equation (\ref{KG}). We 
note that the vector normal form can be simplified near the three 
curves of bifurcations of co-dimension one (see \cite{I95}). The 
truncated normal form (\ref{scalar-fourth-order-equation}) takes 
into account the full problem of bifurcation of co-dimension two.

\section{Proof of existence of center manifold in the system (\ref{center-manifold-Final})--(\ref{resolvent-equation-Final})}

We shall rewrite the system of coupled equations 
(\ref{center-manifold-Final}) and (\ref{resolvent-equation-Final}) 
in the form:
\begin{eqnarray}
\dot{\mbox{\boldmath $\phi$}} &=& {\bf N}(\mbox{\boldmath $\phi$}) + 
\sqrt{\varepsilon} {\bf R}(\mbox{\boldmath $\phi$},\mbox{\boldmath $\psi$}),  \label{basicsyst} \\
\dot{\mbox{\boldmath $\psi$}} & = & \mathcal{L}_{1,0} 
\mbox{\boldmath $\psi$} +\sqrt{\varepsilon } 
\mathbf{F}(\mbox{\boldmath $\phi$},\mbox{\boldmath $\psi$}), 
\label{basicsyst-psi}
\end{eqnarray}%
where $\mbox{\boldmath $\phi$} \in \mathbb{R}^4$ is a function of 
$\zeta_{1} = \sqrt{\varepsilon} \zeta$ and $\mbox{\boldmath $\psi$} 
\in {\cal D}$ is a function of $\zeta$. The vector fields are given 
by
\begin{eqnarray*}
{\bf N}(\mbox{\boldmath $\phi$}) = (\phi _{2},\phi _{3},\phi 
_{4},12\gamma \phi _{3}-12\tau F(\phi _{1}))^{T}, \qquad
\mathbf{F}(\mbox{\boldmath $\phi$},\mbox{\boldmath $\psi$}) = g {\bf 
F}_0(p), 
\end{eqnarray*}
and  
\begin{eqnarray*}
{\bf F}_0 = \left( \begin{array}{cc} 0 \\ -\frac{3}{5} \\ 
\frac{2}{5} p ( 1 - 5 p^2) \end{array} \right), \qquad {\bf 
R}(\mbox{\boldmath $\phi$},\mbox{\boldmath $\psi$}) = \left( 
\begin{array}{cc} 0 \\ -\frac{2}{5} \sqrt{\varepsilon} 
g \\ 0 \\ \frac{12}{\sqrt{\epsilon}} \tilde{g} 
\end{array} \right), 
\end{eqnarray*}
 where $\tilde{g} = g - \gamma \phi_3 + \tau F(\phi_1)$ and $g$ 
is given below (\ref{center-manifold-Final}). We shall rewrite $g$ 
and $\tilde{g}$ in an equivalent form, which is used in our 
analysis. It follows from identities (\ref{difference-operators}), 
representation (\ref{decomposition1})--(\ref{decomposition2}) and 
the scaling (\ref{scaling-transformation}) that 
\begin{eqnarray*}
\delta ^{\pm }U_{3} &=&\phi _{1}\pm \sqrt{\varepsilon }\phi _{2}+\frac{1}{2}%
\varepsilon \phi _{3}\pm \frac{1}{6}\varepsilon ^{3/2}\phi 
_{4}+\varepsilon ^{3/2}\delta ^{\pm }\psi _{3}, \\ \delta _{\pm} U_1 
&=&\phi _{1}(\zeta_1 \pm \sqrt{\varepsilon }) + \varepsilon^{3/2} 
\delta _{\pm} \psi_1, 
\end{eqnarray*}%
such that the equality $\delta^{\pm} U_3 = \delta_{\pm} U_1$ results 
in the identity,
\begin{eqnarray*}
\phi _{3}-2\sqrt{\varepsilon }\psi _{1}+\sqrt{\varepsilon }(\delta 
^{+}+\delta ^{-})\psi _{3} &=&\sqrt{\varepsilon }(\delta_+ \psi_{1} + 
\delta_- \psi _{1}-2\psi _{1})+ \\
&&+\varepsilon ^{-1} \left[ \phi _{1}(\zeta_1 +\sqrt{\varepsilon 
})+\phi _{1}(\zeta_1 - \sqrt{\varepsilon })-2\phi _{1} \right].
\end{eqnarray*}%
As a result, the expression for $g$ below 
(\ref{center-manifold-Final}) can be rewritten as follows:
\begin{eqnarray*}
g = \gamma \sqrt{\varepsilon }(\delta_+ \psi_1 +\delta_- \psi_1 - 
2\psi _{1}) + \gamma \varepsilon ^{-1}(\phi _{1}(\zeta_1 
+\sqrt{\varepsilon})+\phi _{1}(\zeta_1 -\sqrt{\varepsilon })-2\phi 
_{1}) -\tau (\phi _{1}+\varepsilon ^{3/2}\psi _{1}+Q)  
\end{eqnarray*}
and
\begin{eqnarray*}
Q = Q(\phi _{1}(\zeta_1 -\sqrt{\varepsilon })+\varepsilon 
^{3/2}\delta_{-} \psi _{1},\phi _{1}+\varepsilon ^{3/2}\psi 
_{1},\phi _{1}(\zeta_1 +\sqrt{ \varepsilon })+\varepsilon^{3/2} 
\delta_+ \psi _{1}). 
\end{eqnarray*}
Using the system (\ref{basicsyst}), we obtain the following 
representation:
\begin{eqnarray*}
\phi _{1}(\zeta_1 +\sqrt{\varepsilon })+\phi _{1}(\zeta_1 -\sqrt{\varepsilon }%
)-2\phi _{1} &=&\int_{\zeta _{1}}^{\zeta _{1}+\sqrt{\varepsilon 
}}\int_{\tau -\sqrt{\varepsilon }}^{\tau }\ddot{\phi}_{1}(s)dsd\tau \\
&=&\int_{\zeta _{1}}^{\zeta _{1}+\sqrt{\varepsilon }}\int_{\tau -\sqrt{%
\varepsilon }}^{\tau }(\phi _{3}(s) - \frac{2\varepsilon 
}{5}g(s))dsd\tau,
\end{eqnarray*}
such that $g$ is rewritten in an implicit form:
\begin{eqnarray}
\nonumber g &=& \gamma \sqrt{\varepsilon } \left( \delta_+ \psi_1 + 
\delta_- \psi_{1}-2\psi _{1} \right) - \tau (\phi _{1}+\varepsilon ^{3/2}\psi _{1}+Q)  \\
&& + \gamma \varepsilon ^{-1}\int_{\zeta _{1}}^{\zeta _{1}+\sqrt{\varepsilon }%
}\int_{\tau -\sqrt{\varepsilon }}^{\tau }(\phi _{3}(s) - 
\frac{2\varepsilon }{5} g(s))dsd\tau.  \label{def-g} 
\end{eqnarray}
In order to obtain an equivalent expression for $\tilde{g}$, we use 
the system (\ref{basicsyst}) and obtain
\begin{eqnarray*}
-\phi _{3}(\zeta _{1})+\varepsilon ^{-1}\int_{\zeta _{1}}^{\zeta _{1}+\sqrt{%
\varepsilon }}\int_{\tau -\sqrt{\varepsilon }}^{\tau }\phi 
_{3}(s)dsd\tau
=\varepsilon ^{-1}\int_{\zeta _{1}}^{\zeta _{1}+\sqrt{\varepsilon }%
}\int_{\tau -\sqrt{\varepsilon }}^{\tau }\int_{\zeta _{1}}^{s}\phi 
_{4}(\eta )d\eta dsd\tau.
\end{eqnarray*}
It is also clear that 
\begin{eqnarray*}
Q(\phi _{1}(\zeta _{1}-\sqrt{\varepsilon }),\phi _{1}(\zeta 
_{1}),\phi _{1}(\zeta _{1}+\sqrt{\varepsilon }))-Q_{0}(\phi 
_{1}(\zeta _{1}),\phi _{1}(\zeta _{1}),\phi _{1}(\zeta _{1})) \\
= \int_{\zeta _{1}}^{\zeta _{1}+\sqrt{\varepsilon }} \left[ \partial 
_{3}Q(\phi _{1}(\zeta _{1}-\sqrt{\varepsilon }),\phi _{1}(\zeta 
_{1}),\phi _{1}(s))\phi_{2}(s) - \partial _{1}Q(\phi 
_{1}(s-\sqrt{\varepsilon }),\phi _{1}(\zeta _{1}),\phi _{1}(\zeta 
_{1}))\phi _{2}(s-\sqrt{\varepsilon }) \right] ds.
\end{eqnarray*}%
As a result, the expression for $\tilde{g} = g - \gamma \phi_3 + 
\tau F(\phi_1)$ is rewritten in an equivalent form, 
\begin{eqnarray}
\nonumber \tilde{g} &=& \gamma \sqrt{\varepsilon } \left( \delta_+ 
\psi_{1} + \delta_- \psi _{1}-2\psi _{1} \right) - \frac{2}{5} 
\gamma \int_{\zeta _{1}}^{\zeta_{1}+\sqrt{\varepsilon }}\int_{\tau 
-\sqrt{\varepsilon }}^{\tau } g(s)dsd\tau \\
&& - \tau (\varepsilon ^{3/2}\psi _{1}+\widetilde{Q}) + \gamma 
\varepsilon^{-1}\int_{\zeta _{1}}^{\zeta _{1}+\sqrt{\varepsilon }}\int_{\tau -\sqrt{%
\varepsilon }}^{\tau }\int_{\zeta _{1}}^{s}\phi _{4}(\eta )d\eta 
dsd\tau, \label{def-gtilde} 
\end{eqnarray}
where
\begin{eqnarray*}
\widetilde{Q} &=& \int_{\zeta _{1}}^{\zeta _{1}+\sqrt{\varepsilon }} 
\left[ \partial _{3}Q\left[ \phi _{1}(\zeta _{1}-\sqrt{\varepsilon 
}),\phi_{1}(\zeta _{1}),\phi _{1}(s)\right) \phi _{2}(s) \right. \\
&& \left. \phantom{textttexttext} - 
\partial_{1} Q\left( \phi _{1}(s-\sqrt{\varepsilon }),\phi
_{1}(\zeta _{1}),\phi _{1}(\zeta _{1})\right) \phi _{2}(s-\sqrt{\varepsilon }%
)\right] ds \\
&& + Q(\phi_{1}(\zeta_1 -\sqrt{\varepsilon })+\varepsilon^{3/2} 
\delta_- \psi_{1}, \phi _{1}+\varepsilon ^{3/2}\psi _{1},\phi 
_{1}(\zeta_1 +\sqrt{\varepsilon})+\varepsilon ^{3/2}\delta_+ \psi 
_{1}) \\ && - Q(\phi _{1}(\zeta_1 -\sqrt{\varepsilon }),\phi 
_{1},\phi _{1}(\zeta_1 +\sqrt{\varepsilon})).
\end{eqnarray*}%
We observe the fundamental property that for%
\begin{equation}
||\mbox{\boldmath $\phi$}||+||\psi _{1}||\leq M
\end{equation}%
then there is $\tilde{M} > 0$ such that
\begin{equation}
|| g || + \left\|\frac{\tilde{g}}{\sqrt{\varepsilon }}\right\| \leq 
\tilde{M}. \label{bound-g}
\end{equation}
In order to prove Theorem \ref{theorem-center-manifold}, we replace 
the system (\ref{basicsyst-psi}) for $\mbox{\boldmath 
$\psi$}(\zeta)$ by an integral formula. However this is not a simple 
task, since the estimate on the resolvent of $\mathcal{L}_{1,0}$ is 
not nice on its third component. So, it is necessary here to use the 
fact that the nonlinear term in the system
(\ref{basicsyst})--(\ref{basicsyst-psi}) does not depend on 
$\psi_{3}(\zeta,p)$, when the system is rewritten in the new 
formulation (which was precisely the aim of this manipulation). For 
the study of solutions bounded for $\zeta \in \mathbb{R}$, we use 
the two first components of the following implicit formula for 
$\mbox{\boldmath $\psi$}$:
\begin{equation}
\mbox{\boldmath $\psi$}(\zeta) = -\sqrt{\varepsilon} \int_{\zeta}^{\infty} 
e^{(\mathcal{L}_{1,0,+})(\zeta -s)}P_{+}\mathbf{F}(\mbox{\boldmath $\phi$}(s),\mbox{\boldmath $\psi$}(s)) ds 
+\sqrt{\varepsilon }\int_{-\infty }^{\zeta }e^{(\mathcal{L}%
_{1,0,-})(\zeta -s)}P_{-}\mathbf{F}(\mbox{\boldmath 
$\phi$}(s),\mbox{\boldmath $\psi$}(s))ds, \label{integbounded}
\end{equation}
where, by definition
\begin{eqnarray*}
e^{(\mathcal{L}_{1,0,+})t}P_{+}\mathbf{F} &=&\frac{1}{2i\pi 
}\int_{\Gamma _{+}}e^{\Lambda t}(\Lambda 
-\mathcal{L}_{1,0})^{-1}\mathbf{F}d\Lambda, \qquad t<0, \\
e^{(\mathcal{L}_{1,0,-})t}P_{-}\mathbf{F} &=&\frac{1}{2i\pi 
}\int_{\Gamma _{-}}e^{\Lambda t}(\Lambda 
-\mathcal{L}_{1,0})^{-1}\mathbf{F}d\Lambda , \qquad t>0,
\end{eqnarray*}%
The curves $\Gamma _{\pm }$ are defined by 
\begin{eqnarray*}
\Gamma _{+} &=&C_{+}\cup \overline{C_{+}}\cup L_{+}, \\
\Gamma _{-} &=&C_{-}\cup \overline{C_{-}}\cup L_{-},
\end{eqnarray*}
where $C_{\pm }$ are curves on the complex plane defined by
\begin{eqnarray*}
C_{\pm } = \left\{ \Lambda = \pm x+iy : \; y = c \cosh (\alpha x), 
\; \alpha >1, \; c > \sqrt{8}, \; |x|\geq \xi >0 \right\},
\end{eqnarray*}%
and
\begin{equation}
L_{\pm} = \left\{ \Lambda =\pm \xi +iy : \; |y|\leq c \cosh (\alpha 
\xi ) \right\}.
\end{equation}%
The curve $\Gamma _{-}$ is oriented with increasing ${\rm 
Im}(\Lambda)$, while the curve $\Gamma _+$ is oriented with 
decreasing ${\rm Im}(\Lambda)$. It is clear from the proof of Lemma 
\ref{lemma-eigenvalues} that both curves $\Gamma _{\pm }$ lie in the 
resolvent set of $\mathcal{L}_{1,0}$ and that on these curves we 
have (as on the imaginary axis for $\Lambda =ik$, $|k|$ large 
enough),
\begin{equation}
|D_{1}|\geq c_{1}|\Lambda |^{2},
\end{equation}
where $D_1(\Lambda) = 2 ( \cosh \Lambda - 1) - \Lambda^2$. We shall 
complete the system (\ref{basicsyst-psi}) with an explicit formula 
for $\psi_{3}(\zeta,p)$. In order to study integrals in the integral 
equation (\ref{integbounded}), we compute the expression 
\begin{equation}
\mathbf{R}(\Lambda )=(\Lambda -\mathcal{L}_{1,0})^{-1} {\bf 
F}_{0}(p).
\end{equation}
Using solutions of the resolvent equation
(\ref{solution-resolvent-equation-0})--(\ref{solution-resolvent-equation}) 
for $\mathbf{F} = {\bf F}_0(p)$, we obtain the following 
expressions:
\begin{equation}
\widetilde{F}(\zeta,\Lambda) = -\frac{3}{5}-\frac{4}{5\Lambda 
^{2}}(14+\cosh \Lambda ) + \frac{24}{\Lambda ^{4}} (\cosh \Lambda - 
1),
\end{equation}
and
\begin{eqnarray*}
\int_{0}^{p}\frac{2}{5}s(1-5s^{2})e^{\Lambda (p-s)}ds = -\frac{2p(1-5p^{2})%
}{5\Lambda }+\frac{2(15p^{2}-1)}{5\Lambda ^{2}}+\frac{2e^{\Lambda p}}{%
5\Lambda ^{2}}+\frac{12p}{\Lambda ^{3}}+\frac{12}{\Lambda 
^{4}}(1-e^{\Lambda p}).
\end{eqnarray*}
As a result, the first two components of $\mathbf{R}(\Lambda)$ 
satisfy the following bounds for $\Lambda \in \Gamma _{\pm }$:
\begin{eqnarray*}
\left\vert R_{1}(\Lambda )-\frac{3}{5\Lambda ^{2}} \right\vert &\leq 
& \frac{C_1}{|\Lambda|^{4-1/\alpha }} \\
\left\vert R_{2}(\Lambda )-\frac{3}{5\Lambda } \right\vert &\leq & 
\frac{C_2}{|\Lambda |^{3-1/\alpha }}
\end{eqnarray*}
for some constants $C_{1,2} > 0$. Since the following integrals 
cancel for $t<0$ (below we consider integrals on $\Gamma _{+},$ 
analogous results hold on $\Gamma _{-}$ for $t>0$),
\begin{eqnarray*}
\frac{1}{2i\pi }\int_{\Gamma _{+}}\frac{e^{\Lambda t}}{\Lambda }d\Lambda =%
\frac{1}{2i\pi }\int_{\Gamma ^{+}}\frac{e^{\Lambda t}}{\Lambda 
^{2}}d\Lambda =0,
\end{eqnarray*}
we obtain 
\begin{equation}
\left\vert \frac{1}{2i\pi }\int_{\Gamma _{+}}e^{\Lambda 
t}R_{j}(\Lambda )d\Lambda \right\vert \leq \tilde{C}_j \frac{e^{-\xi 
(\alpha -t)}}{\alpha -t}
\end{equation}
and
\begin{equation}
\left\vert \frac{d}{dt}\left( \frac{1}{2i\pi }\int_{\Gamma 
_{+}}e^{\Lambda t}R_{j}(\Lambda )d\Lambda \right) \right\vert \leq 
\tilde{C}'_j \frac{e^{-\xi (\alpha -1-t)}}{\alpha -1-t},
\end{equation}
where $j=1,2$ and $t < \min(0,\alpha-1)$ (note that $\alpha 
> 1$). As a consequence, for $\alpha >1$, the two first components of the 
integral $\frac{1}{2i\pi }\int_{\Gamma _{+}}e^{\Lambda 
t}\mathbf{R}(\Lambda )d\Lambda $ and $\frac{1}{2i\pi }\int_{\Gamma 
_{-}}e^{\Lambda t}\mathbf{R}(\Lambda )d\Lambda $ are differentiable 
real functions respectively for $t\leq 0$ and for $t\geq 0$. We note 
that no similar estimates can be obtained for the third component of 
$\mathbf{R}(\Lambda)$ for $t<0$. Now using the identity
\begin{equation}
\Lambda (\Lambda -\mathcal{L}_{1,0})^{-1}\mathbf{F}=(\Lambda -\mathcal{L}%
_{1,0})^{-1}\mathcal{L}_{1,0}\mathbf{F}+\mathbf{F}
\end{equation}%
one can show that a bounded solution $(\psi_{1},\psi_{2})\in 
C_{b}(\mathbb{R}, \mathbb{R}^{2})$ of the two first components of 
the integral equation (\ref{integbounded}) is completed by the 
following explicit formula for $\psi _{3}(\zeta,p)$
\begin{equation}
\psi _{3}(\zeta ,p)=\psi _{1}(\zeta +p) - \frac{\sqrt{\varepsilon }}{10}%
\int_{\zeta -p}^{\zeta +p}(\zeta +p-\tau )\left( 1-\frac{5}{4}(\zeta 
+p-\tau )^{2}\right) g\left( \frac{\zeta +p+\tau }{2} \right)d\tau,  
\label{psi_3}
\end{equation}%
such that $\psi _{3}(\zeta ,0)=\psi _{1}(\zeta)$. Therefore, for 
bounded solutions on $\zeta \in \mathbb{R}$, the system 
(\ref{basicsyst-psi}) with the relation (\ref{psi_3}) reduces to a 
two-dimensional system. The argument of \cite[p.145]{VI92} can then 
be used.

Coming back to the system (\ref{basicsyst})--(\ref{basicsyst-psi}), 
the implicit linear equation (\ref{def-g}) for $g$ can be solved for 
small enough $\varepsilon$ and the system 
(\ref{basicsyst})--(\ref{basicsyst-psi}) is hence equivalent to the 
coupled system 
(\ref{center-manifold-Final})--(\ref{resolvent-equation-Final}). 
From the above study, the two first components of equation 
(\ref{integbounded}), completed by (\ref{psi_3}) have a meaning for 
$\mbox{\boldmath $\psi$}$ in $C_{b}^{0}(\mathbb{R},\mathcal{D})$, 
for any known function $\mbox{\boldmath $\phi$} \in 
C_{b}^{0}(\mathbb{R},\mathbb{R}^{4})$. Moreover, for $\varepsilon $ 
small enough, and thanks to (\ref{bound-g}), the implicit function 
theorem applies, and one can find a unique solution of the integral 
equation (\ref{integbounded}) for $\mbox{\boldmath $\psi$} \in
C_{b}^{0}(\mathbb{R},\mathcal{D})$, which is bounded by a constant 
of the order of ${\rm O}(\sqrt{\varepsilon })$ under the condition 
that 
\begin{equation}
||\mbox{\boldmath $\phi$} 
||_{C_{b}^{0}(\mathbb{R},\mathbb{R}^{4})}\leq M,
\end{equation}%
where $M$ is an arbitrarily fixed number. Moreover, because of the 
reversibility of the system (\ref{basicsyst})--(\ref{basicsyst-psi}) 
and oddness property of the function $Q(v,u,w)$, if $\mbox{\boldmath 
$\phi$}(\zeta_1)$ is anti-reversible with respect to ${\cal S}_0$, 
then the unique solution $\mbox{\boldmath $\psi$}(\zeta)$ is also 
anti-reversible with respect to ${\cal S}$. 

It remains to prove the last sentence of Theorem 
\ref{theorem-center-manifold}. When $\mbox{\boldmath 
$\phi$}(\zeta_1)$ is a periodic function, there exists a unique 
periodic solution $\mbox{\boldmath $\psi$}(\zeta)$ of the system 
(\ref{basicsyst-psi}), since the system 
(\ref{basicsyst})--(\ref{basicsyst-psi}) is invariant under a shift 
of $\zeta$ by $T/\sqrt{\varepsilon }.$ When $\mbox{\boldmath 
$\phi$}(\zeta_1)$ is an asymptotically periodic function at 
infinity, the solution $\mbox{\boldmath $\psi$}(\zeta)$ is then 
asymptotic to the corresponding periodic solution of the system 
(\ref{basicsyst-psi}). This follows from the fact that the formula 
(\ref{psi_3}) for the third component of $\mbox{\boldmath 
$\psi$}(\zeta)$ can be considered in the limit $\zeta \rightarrow 
\pm \infty $ and from the two-dimensional equation 
(\ref{integbounded}) for the two first components of 
$\mbox{\boldmath $\psi$}(\zeta)$, for which the solution depends 
continuously on $\mbox{\boldmath $\phi$}(\zeta_1)$, as for an 
ordinary differential equation.

\section{Stokes constant for the fourth-order ODE (\ref{ODE})}
                  
The Stokes constants are used in rigorous analysis of non-existence 
of homoclinic and heteroclinic orbits in singularly perturbed 
systems of differential and difference equations \cite{T00}. If the 
Stokes constant is non-zero, the heteroclinic orbit of the 
unperturbed problem does not persist beyond all orders of the 
perturbation series expansion. We will develop here computations of 
the Stokes constant for the fourth-order equation (\ref{ODE}) in the 
limit $\sigma \to \infty$, in the context of the $\phi^4$ model. We 
will show that the Stokes constant is non-zero, which implies that 
the graph of $K(\sigma)$ on Figure 3(a) remains non-zero in the 
limit $\sigma \to \infty$. 
                           
Let us introduce a small parameter $\epsilon$, such that $\epsilon = 
\sigma^{-1}$, and apply rescaling of $t$ as $x = \sqrt{\epsilon} t$. 
The fourth-order equation (\ref{ODE}) with $F = \phi ( 1 - \phi^2)$ 
takes the form of the singularly perturbed ODE:
\begin{equation}
\label{ODE-appendix} \phi'' + \phi - \phi^3 + \epsilon^2 \phi^{({\rm 
iv})} = 0.
\end{equation}
At $\epsilon = 0$, the ODE (\ref{ODE-appendix}) corresponds to the 
stationary problem for the continuous $\phi^4$ model with the 
heteroclinic orbit in the form $\phi_0(x) = \tanh(x/\sqrt{2})$. 
There exists a formal solution to the perturbed ODE 
(\ref{ODE-appendix}) in the form of the regular perturbation series,
\begin{equation}  
\label{formal-solution-appendix} \hat{\phi}(x) = \tanh 
\frac{x}{\sqrt{2}} + \sum_{k \geq 1} \epsilon^{2k} \phi_k(x),
\end{equation}
where an odd function $\phi_k(x)$ is uniquely defined from the 
solution of the linear inhomogeneous problem:
\begin{equation}
\left( - \partial_x^2 + 2 - 3 {\rm sech}^2 \frac{x}{\sqrt{2}} 
\right) \phi_k = H_k(\phi_{k-1},\phi_{k-2},...,\phi_0), \qquad k 
\geq 1,
\end{equation}
where $H_k$ are correction terms which are odd exponentially 
decaying functions of $x \in \mathbb{R}$. Each term of $\phi_k(x)$ 
is real analytic function of $x \in \mathbb{R}$ and is extended 
meromorphically to the complex plane, with pole singularities at the 
points $x = \frac{\pi (1 + 2 n) i}{\sqrt{2}}$, $n \in \mathbb{Z}$. 
The Stokes constants are introduced after the rescaling of the 
perturbed ODE (\ref{ODE-appendix}) near the pole singularities and 
replacing the formal solution (\ref{formal-solution-appendix}) with 
the inverse power series. We rescale the solution as follows: 
$$
x = \frac{\pi i}{\sqrt{2}} + \epsilon z, \qquad \phi(x) = 
\frac{\psi(z)}{\epsilon},
$$
where $\psi(z)$ satisfies the regularly perturbed ODE:
\begin{equation}                                        
\label{ODE-regular-appendix} \psi^{({\rm iv})} + \psi'' - \psi^3 + 
\epsilon^2 \psi = 0.
\end{equation}          
We are looking at the formal solution of the ODE
(\ref{ODE-regular-appendix}) with the perturbation series 
\begin{equation}  
\label{inner-solution-appendix} \hat{\psi}(z) = \hat{\psi}_0(z) + 
\sum_{k \geq 1} \epsilon^{k} \hat{\psi}_k(z),
\end{equation}                             
where each term $\hat{\psi}_k(z)$ is the inverse power series of 
$z$, starting with the zero-order solution:
\begin{equation}
\label{zero-solution-appendix} \hat{\psi}_0(z) = \sum_{m \geq 1} 
\frac{a_m}{z^m}.
\end{equation}
Matching conditions between (\ref{formal-solution-appendix}), 
(\ref{inner-solution-appendix}), and (\ref{zero-solution-appendix}) 
imply that $a_1 = \sqrt{2}$ and $a_2 = 0$. The main result of the 
beyond-all-order perturbation theory is that the formal series 
(\ref{formal-solution-appendix}) are (\ref{inner-solution-appendix}) 
diverges if the inverse power series (\ref{zero-solution-appendix}) 
diverges. Divergence of the formal series implies in its turn that 
the heteroclinic orbit  $\phi_0(x) = \tanh(x/\sqrt{2})$ is destroyed 
by the singular perturbation of the ODE (\ref{ODE-appendix}). 
Divergence of (\ref{zero-solution-appendix}) is defined by the 
asymptotical behavior of the coefficients $a_m$ as $m \to \infty$ 
(see \cite{TP05} for details). 

In order to show that the inverse power series 
(\ref{zero-solution-appendix}) diverges, we find the recurrence 
relation between coefficients in the set $\{ a_m \}_{m=1}^{\infty}$:
\begin{equation} 
\label{recurrence-equation-appendix} m (m+1)(m+2)(m+3) a_m + 
(m+2)(m+3) a_{m+2} - \sum_{l + k \leq m+3} a_l a_k a_{m+4-l-k} = 0, 
\end{equation}                           
where $m \geq 1$, $a_1 = \sqrt{2}$ and $a_2$ is arbitrary. Due to 
the translational invariance, we can always fix $a_2 = 0$. Setting 
$a_m = (-1)^n (2 n)! b_n$ for $m = 2 n + 1$, $n \geq 0$ and $a_m = 
0$ for $m = 2n$, $n \geq 1$, we reduce the recurrence equation 
(\ref{recurrence-equation-appendix}) to the diagonal form:
\begin{equation}
b_{n+1} = b_n + \sum_{l + k \leq n+1} \frac{(2l)! (2k)! 
(2n+2-2l-2k)!}{(2n+4)!} b_l b_k b_{n+1-l-k}, \qquad n \geq 0,
\end{equation}                                               
where $b_0 = \sqrt{2}$. Since 
$$
1 - \frac{6}{(2n+4)(2n+3)} > 0, \qquad n \geq 0,
$$
we prove that the sequence $\{ b_n \}_{n = 0}^{\infty}$ is 
increasing, such that $b_{n+1} \geq b_n > b_0 = \sqrt{2}$. 
Therefore, the sequence $\{ b_n \}_{n = 0}^{\infty}$ is bounded from 
below by a positive constant, so that the Stokes constant in the 
beyond-all-order theory is strictly positive.

Similar computations can be developed for the fourth-order equation 
(\ref{ODE}) with $F = \sin(\phi)$. They are left for a reader's 
exercise.

\section{The normal form for the inverse method of \cite{FZK99}}

It was shown in \cite{FZK99} that the nonlinearity of the discrete 
Klein--Gordon equation (\ref{Klein-Gordon}) can be chosen in such 
way that it yields an exact travelling kink solution $u_n(t) = 
\phi(z)$, $z = n-ct$ for a particular value of velocity $c = s$. In 
application to the discrete $\phi^4$ lattice, one can look for the 
exact solution in the form $\phi(z) = \tanh(\mu z)$, where $\mu$ is 
arbitrary parameter, and derive the explicit form of the 
nonlinearity function $f = f(u_n)$ parameterized by $s$ and $\mu$. 
(See eq. (27) and refs. [13,14] in \cite{FZK99}.) Since we are using 
a particular form of the nonlinearity $f(u_{n-1},u_n,u_{n+1})$ in 
the starting equation (\ref{Klein-Gordon}), we transform the 
nonlinearity from \cite{FZK99} to the form:
\begin{equation} 
\label{nonlinearity-FZK} f(u_n) = u_n ( 1 - u_n^2) \frac{1+ \alpha 
u_n^2}{1 - \beta u_n^2},
\end{equation}
where 
$$
\alpha = \frac{s^2 \mu^2 \tanh^2 \mu}{\tanh^2 \mu - s^2 \mu^2}, 
\qquad \beta = \tanh^2 \mu,
$$
and $h^2$ is related to the parameters $\mu$ and $s$ by means of the 
relation:
\begin{equation}
\label{transcendental-appendix} h^2 = 2 \left( \tanh^2 \mu - s^2 
\mu^2 \right).
\end{equation}  
We note that the nonlinearity (\ref{nonlinearity-FZK}) violates the 
assumption (\ref{nonlinearity-continuous}) of our paper, since it 
depends implicitly on the lattice parameter $h^2$. 

When $s = 0$ ($\alpha = 0$, $\beta = \frac{1}{2} h^2$), the 
nonlinearity (\ref{nonlinearity-FZK}) represents one of the four 
exceptional nonlinearities with translationally invariant stationary 
kinks (see \cite{BOP05}). We consider here the limit $s \to 1$ from 
the point of normal form analysis. We will show that existence of 
the exact travelling kink solutions $\phi(z) = \tanh(\mu z)$ in the 
differential advance-delay equation (\ref{advance-delay}) is 
preserved within the reduction to the scalar fourth-order equation 
(\ref{scalar-fourth-order}). Let 
$$
s^2 = 1 + \epsilon \gamma_s, \quad c^2 = 1 + \epsilon \gamma, \qquad 
h^2 = \epsilon^2 \tau.
$$
The parameter $\mu$ is defined from the transcendental equation 
(\ref{transcendental-appendix}) in the form, 
$$
\mu = \sqrt{\epsilon} \mu_s + {\rm O}(\epsilon),
$$
where $\mu_s$ is a real positive root of bi-quadratic equation:
$$
\frac{4}{3} \mu_s^4 + 2 \gamma_s \mu_s^2 + \tau = 0.
$$
Two real positive roots $\mu_s$ exist in the domain $0 < \tau < 
\frac{3}{4} \gamma_s^2$ for $\gamma_s < 0$ and no real positive 
roots exist beyond this domain for $\tau > 0$. The nonlinearity 
(\ref{nonlinearity-FZK}) is reduced in the limit of small $\epsilon$ 
to the form,
$$
F(u) = u(1-u^2)( 1 + \frac{2 \mu_s^4}{\tau} u^2) + {\rm 
O}(\epsilon),
$$
such that the normal form (\ref{scalar-fourth-order}) becomes now
\begin{equation}              
\label{appendix-normal-form} \frac{1}{12} \phi^{({\rm iv})} - \gamma 
\phi'' + \tau \phi (1 - \phi^2) + 2 \mu_s^4 \phi^3 ( 1 - \phi^2) = 
0. 
\end{equation}                               
It is easy to verify that the scalar fourth-order equation 
(\ref{appendix-normal-form}) has the exact heteroclinic orbit 
$\phi(z) = \tanh( \mu_s z)$ for $\gamma = \gamma_s$ and $\tau > 0$. 
Since there are two roots for $\mu_s$, two heteroclinic orbits 
$\phi(z)$ exist. These two orbits disappear at the saddle-node 
bifurcation as $\tau \to \frac{3}{4} \gamma_s^2$. Within the 
numerical method of our paper, the heteroclinic orbits can be 
detected as isolated zeros of the split function $K(\sigma)$. We 
conclude that modifications of assumptions on the nonlinearity of 
the discrete Klein--Gordon equation (\ref{Klein-Gordon}) result in 
modifications of the scalar fourth-order normal form 
(\ref{scalar-fourth-order}), which may hence admit isolated 
heteroclinic orbits.

\end{document}